\documentstyle[12pt]{article}
\begin{document}
\author{ S.V. Ludkovsky}
\title{Absolute non-Archimedean polyhedral expansions
of ultrauniform spaces.}
\date{17 May 2000}
\maketitle
\par Theoretical Department,
Institute of General Physics, Str. Vavilov 38, 
Moscow, 117942, Russia
\section{Introduction.}
\par To expansions of topological and uniform spaces
as inverse limits several works were devoted 
\cite{coris,fed,freu,isb1,isb2,isb3,isb4,isb5,koz,pas}.
Irreducible polyhedral expansions 
over the field $\bf R$ were very important in \cite{freu,isb3,
koz}. In the case of locally compact groups expansions
into inverse limits of Lie groups were investigated in
\cite{monzi,pont}. Such inverse mapping systems expansions
may be used for the proof of isomorphisms of definite
locally convex spaces over the fields $\bf R$ or $\bf C$
\cite{lu1,lu2,lu3}. In the case of ultrauniform hereditarily
disconnected spaces polyhedral expansions over $\bf R$ do not take 
into account many features of such spaces.
\par This work is devoted to the investigation of the problem about
inverse mapping systems 
expansions of ultrauniform spaces using polyhedra over 
non-Archimedean locally compact fields. In \S 2 preliminary results,
definitions and notations are given. In \S 3 Theorems about expansions
of complete ultrametric and ultrauniform spaces are proved. 
In \S 4 absolute polyhedral expansions and inverse mapping systems
of expansions
for non-complete spaces are investigated. There are elucidated cases,
when in inverse mapping systems
of expansions polyhedra are finite-dimensional
over non-Archimedean fields $\bf L$. In particular locally compact 
ultrauniform groups are considered. The main results are Theorems 3.18
and in \S 4.
\par The considered here topological spaces $X$ are hereditarily disconnected
and zero-dimensional, but in general they are not necessarily strongly 
zero-dimensional, that is this class has representatives with
covering dimensions from $0$ to $\infty $ 
(see Chapters 6 and 7 in \cite{eng}). This article also contains 
results about a relation of $dim (X)$ and dimensions of polyhedra over
$\bf L$. With the help of polyhedral expansions it is shown in \S 5 
that there exists a natural correspondence between ultrauniform 
$X$ and uniform spaces $Y$, such that there exists a continuous 
quotient mapping $\theta : X\to Y$. That is for each $X$ there exist
$\theta $ and $Y$, vice versa for each $Y$ there exist $X$ and $\theta $
such that $d(X)=d(Y)$ and $w(X)=w(Y)$. 
Then $X$ is dense in itself if and only
if $Y$ is such. Moreover, there are natural compactifications
$cX$ and $c_{\theta }Y$ such that $\theta $ has an extension on 
$cX$ and $\theta (cX)=c_{\theta }Y$. On the other hand, it gives a
relation between dimensions $dim_{\bf L}P$ of polyhedras $P$ over $\bf L$
in an expansion and $dim [\theta (P)]$,
also between $\max_{P_m} dim_{\bf L}P_m$ and $dim [\theta (X)]$.
Then an important particular case of polyhedral expansions
of manifolds modelled on non-Archimedean locally convex spaces 
is investigated. Few results from \S 3 and \S 4 were announced in
\cite{lu5}.
\section{Notations and preliminary results.}
\par {\bf 2.1.} Let us recall that by an ultrametric space $(X, \rho )$ 
is implied a set $X$ with a metric $\rho $ such that
it satisfies the ultrametric inequality:
$\rho (x,y)\le \max (\rho (x,z); \rho (z,y) )$ for each $x,$ $y$
and $z\in X$. A uniform space $X$ with ultrauniformity $U$ is called
an ultrauniform space $(X, U)$ such that $U$ satisfies the following
condition: $|x-z|<V'$, if $|y-z|<V'$ and $|x-y|<V$, where
$V\subset V'\in U$, $x,$ $y$ and $z\in X$ \cite{eng,roo}.
If in $X$ a family of pseudoultrametrics $\bf P$ is given and it satisfies
conditions $(UP1,UP2)$, then it induces an ultrauniformity $U$
due to Proposition 8.1.18 \cite{eng}.
\par Let $\bf L$ be a non-Archimedean field \cite{roo}. We say that $X$ is 
a $\bf L$-Tychonoff space, if $X$ is a $T_1$-space and for each
$F= \bar F \subset X$ with $x\notin F$ there exists a continuous function
$f: X\to B({\bf L},0,1)$ such that $f(x)=0$, $f(F)=\{ 1\} $, where
$B(X,y,r):= \{ z\in X: \rho (y,z)\le r \} $ for $y\in X$ and $r\ge 0$.
From $ind({\bf L})=0$ it follows that the small inductive dimension
is $ind(X)=0$ (see \S 6.2 and ch. 7
in \cite{eng}).
Since the norm $|*|_{\bf L}: {\bf L}\to \tilde \Gamma _{\bf L}$ 
is continuous, then $X$ is the Tychonoff space, where 
$\tilde \Gamma _{\bf L}:=\{ |x|_{\bf L}: x\in {\bf L} \} \subset 
[0, \infty )$. Vice versa if $X$ is a Tychonoff space with $ind(X)=0$,
then it is also $\bf L$-Tychonoff, since there exists a clopen
(closed and open at the same time) neighbourhood $W\ni x$ with
$W\cap F=\emptyset $ and as the locally constant function $f$ may 
be taken with $f(x)=0$ and $f(X\setminus W)= \{ 1\} $. 
\par Let us consider spaces $C(X,{\bf L}):= \{ f: X\to {\bf L}|
f \mbox{ is continuous } \} $ and $C^*(X,{\bf L}):=\{ f\in C(X,{\bf L}):
|f(X)|_{\bf L} \mbox{ is bounded in } {\bf R} \} $, then for each finite
family $\{ f_1,...,f_m\} \subset C(X,{\bf L})$ 
(or $C^*(X,{\bf L})$) the following
pseudoultrametric is defined: $\rho _{f_1,...,f_m}(x,y):=\max (
|f_j(x)-f_j(y)|_{\bf L}: j=1,...,m \} $. Families $\bf P$ or $\bf P^*$
of such $\rho _{f_1,...,f_m}$ induce ultrauniformities $\sf C$ or $\sf C^*$
respectively and the initial topology on $X$. If a sequence 
$\{ V_j: j=0,1,... \}
\subset U$ is such that $V_0=X^2$, $pV_{j+1}\subset V_j$ for $j=1,2,...,$
where $p$ is a prime number, then there exists a pseudoultrametric
$\rho $ such that
$\rho (x,y):=0$ for $(x,y)\in \bigcap_{j=0}^{\infty }V_j$, $\rho (x,y)=
p^{-j}$ for $(x,y)\in V_j\setminus V_{j+1}$, so $V_i\subset
\{ (x,y): \rho (x,y)\le p^{-i} \} \subset V_{i-1}$. Indeed, from $(x,y)
\in V_i\setminus V_{i+1}$ and $(y,z)\in V_j\setminus V_{j+1}$ for $j\ge i$ 
it follows $(x,z)\in V_i$ and $\rho (x,z)\le p^{-i}=\rho (x,y)$.
Therefore, ultrauniform spaces may be equivalently characterised 
by $U$ or $\bf P$ (see \S 8.1.11 and \S 8.1.14 in \cite{eng} ).
\par Henceforth, locally compact non-discrete non-Archimedean
infinite fields $\bf L$ are considered. If the characteristic
$char({\bf L})=0$, then due to 
\cite{wei} for each such $\bf L$ there exists a prime number $p$
such that $\bf L$ is a finite algebraic extension of the field
$\bf Q_p$ of $p$-adic numbers. If $char({\bf L})=p>0$, then $\bf L$
is isomorphic with the field $\bf F_p(\theta )$ of the formal power series
by the indeterminate $\theta $, where $p$ is the prime number,
each $z\in \bf L$ has the form $z=\sum_{j=k(z)}^{+ \infty }
a_j\theta ^j$ with $a_j\in \bf F_p$ for each $j$, 
$|z|_{\bf L}=p^{-k(z)}$
up to the equivalence of the non-Archimedean valuation,
$k(z)\in \bf Z$, $\bf F_p$ is the finite field consisting of
$p$ elements, ${\bf F_p}\hookrightarrow \bf L$ is the natural 
embedding.
\par For an ordinal $A$ with its cardinality
$m=card(A)$ by $c_0({\bf L},A)$ it is denoted the following Banach space
with vectors $x=(x_a: a\in A, x_a\in {\bf L})$ of a finite norm
$\| x\|:=\sup_{a\in A}|x_a|_{\bf L}$ and such that for each $b>0$
a set $\{ a\in A: |x_a|_{\bf L}\ge b \} $ is finite. It has the orthonormal
in the non-Archimedean sense basis 
$\{ e_j:=(\delta _{j,a}: a\in A): j\in A\} $,
where $\delta _{j,a}=1$ for $j=a$ and $\delta _{j,a}=0$ for
each $j\ne a$ \cite{roo}. If $card(A_1)=card(A_2)$ then $c_0({\bf L},A_1)$
is isomorphic with $c_0({\bf L},A_2)$. In particular $c_0({\bf L},n)=\bf L^n$
for $n\in \bf N$. Then $card(A)$ is called the dimension of $c_0({\bf L},A)$
and it is denoted by $card(A)=dim_{\bf L}c_0({\bf L},A)$.
\par Let $\sf U$ be a uniform covering of $(X,{\bf P})$
(see \S 8.1 in \cite{eng}), then the collection
$\{ St(U,{\sf U}): U\in {\sf U} \}  =: \sf \mbox{ }^*U$
is called the star of $\sf U$. If $\sf \mbox{ }^*U$
is a refinement of a uniform covering $\sf V$, then $\sf U$ is called 
the star refinement of $\sf V$.
\par {\bf 2.2. Lemma.} {\it Let $(X, \rho )$ be an ultrametric space,
then there exists an ultrametric $\rho '$ equivalent to $\rho $ such that
$\rho '(X,X)\subset \tilde \Gamma _{\bf L}$.}
\par {\bf Proof.} Let $\rho '(x,y):=\sup_{b\in \tilde \Gamma _{\bf L},
b\le \rho (x,y)}b$, 
where $x$ and $y\in X$, either $\bf L\supset Q_p$
or ${\bf L}=\bf F_p(\theta )$.
Then $\rho '$ is the ultrametric such that
$\rho '(x,y)\le \rho (x,y)\le p\times \rho '(x,y)$ for 
each $x$ and $y\in X$.
\par {\bf 2.3. Lemma.} {\it Let $(X,{\bf P})$ be an ultrauniform space,
then there exists a family $\bf P'$ such that $\rho '(X,X)\subset 
\tilde \Gamma _{\bf L}$ for each $\rho '\in \bf P'$; $(X,{\bf P})$ and
$(X, {\bf P'})$ are uniformly isomorphic, the completeness of one of them 
is equivalent to the comleteness of another.}
\par {\bf Proof.} In view of Lemma 2.2 for each $\rho \in \bf P$ there exists
an equivalent pseudoultrametric $\rho '$. They form a family $\bf P'$.
Evidently, the identity mapping $id: (X, {\bf P})\to (X,{\bf P'})$
is the uniform isomorphism. The last statement follows from
8.3.20 \cite{eng}.
\par {\bf 2.4. Theorem.} {\it For each ultrametric space $(X,\rho )$
there exist an embedding $f: X\to B(c_0({\bf L},A_X),0,1)$ and
an uniformly continuous embedding into $c_0({\bf L},A_X)$, where
$card (A_X)=w(X)$, $w(X)$ is the topological weight of $X$.}
\par {\bf Proof.} In view of Theorem 7.3.15 \cite{eng}
there exists an embedding of $X$ into the Baire space $B(m)$, where
$m=w(X)\ge \aleph _0$. In the case $w(X)<\aleph _0$ this statement
is evident, since $X$ is finite. In view of Lemma 2.2 we choose
in $B(m)$ an ultrametric $\rho $ equivalent to the initial one with 
values in $\tilde \Gamma _{\bf L}$ such that $\rho (\{ x_i\},
\{ y_i\} )=p^{-k}$, if $x_k\ne y_k$ and $x_i=y_i$ for $i<k$,
$\rho (\{ x_i\} ,\{ y_i\} )=0$, if $x_i=y_i$ for all $i$, where 
$ \{ x_i\} \in B(m)$, $i\in \bf N$, $x_i\in D(m)$, $D(m)$ 
denotes the discrete space of cardinality $m$. Let 
$A={\bf N}\times C$, $card (C)=m$, $\{ e_{i,a}| (i,a)\in A\}$ 
be the orthonormal basis in $c_0({\bf L},A)$. For each $\{ x_i\} 
\in B(m)$ we have $x_i\in D(m)$ and we can take $x_i=e_{i,a}$ for 
suitable $a=a(i)$, since $D(m)$ is isomorphic with $\{ e_{i,a}: a\in C \} $.
Let $f(\{ x_i\}):=\sum_{i\in {\bf N}, a=a(i)}p^ie_{i,a}$, consequently,
$\| f( \{ x_i\} )-f(\{ y_i\} )\|_{c_0({\bf L},A)}=\rho (\{ x_i\} ,\{ y_i\} )$.
\par The last statement of the Theorem follows
from the isometrical embedding of $(X,\rho )$
into the corresponding free Banach space, which is isomorphic with
$c_0({\bf L},A)$ (see Theorem 5 in \cite{lu3} and Theorems 5.13 and
5.16 in \cite{roo}), since each Banach space over a non-Archimedean 
spherically complete field $\bf L$ is isomorphic with $c_0({\bf L},A)$
with the corresponding $A$. Instead of the free Banach space generated by
$(X,\rho ')$ (see also Lemma 2.2), 
it can be used an embedding into the topologically dual space
${\tilde C}^0_b(X,{\bf L})^*$ to a space ${\tilde C}^0_b(X,{\bf L})$,
which is the space of bounded continuous 
Lipschitz functions $f: X\to \bf L$
with $\| f\| :=\max \{ \sup_{x\in X} |f(x)|, \sup_{x\ne y\in X}
(|f(x)-f(y)|/\rho (x,y)) \} <\infty .$ 
For each $S\subset X$ there exists a function $f\in 
{\tilde C}^0_b(X,{\bf L})$ such that $\rho '(x,S)/[1+\rho '(x,S)]\le |f(x)|
\le p\rho '(x,S)/[1+\rho '(x,S)]$, where $\rho '(x,S):=\inf_{y\in S}
\rho '(x,y)$.
To each $x\in X$ there corresponds
a continuous functional $G_x$ such that $G_x(f)=f(x)$, so 
$G_x(f)-G_y(f)=f(x)-f(y)$. In view of the $\bf L$-regularity
of $X$ we have $\| G_x-G_y \| =\rho (x,y)$.
\par {\bf 2.5. Corollary.} {\it Each ultrauniform space $(X,{\bf P})$ 
has a topological embedding into $\prod_{\rho \in \bf P}
B(c_0({\bf L},A_{\rho }),0,1)$ and an uniform embedding into
$\prod_{\rho \in \bf P}c_0({\bf L},A_{\rho })$
with $card(A_{\rho })=w(X_{\rho })$.}
\par {\bf Proof.} For each ultrametric $\rho \in \bf P$ 
of an ultrauniform space
$(X,{\bf P})$ there exists the equivalence relation $R_{\rho }$
such that $xR_{\rho }y$ if and only if $\rho (x,y)=0$. Then there exists
the quotient mapping $g_{\rho }: X\to X_{\rho }:=
(X/R_{\rho })^{~}$, where $X_{\rho }$ is 
the ultrametric space with the ultrametric also denoted by $\rho $,
$\tilde X$ denotes the completion of $X$. Then $X$ has the uniform
embedding into the limit of the inverse spectra 
$\lim_{\rho }{\tilde X}_{\rho }=\tilde X$. 
\par {\bf 2.6. Corollary.} {\it For each $j\in \bf N$ and $f(X)
\subset c_0({\bf L},A_X)$ from Theorem 2.4 there are coverings $U_j$ of 
the space $f(X)$ by disjoint clopen balls $B_l$ with diameters not greater
than $p^{-j}$ and with $\inf_{l\ne k}dist(B_l,B_k)>0$.}
\par {\bf Proof} follows from the consideration of the covering
of the Banach space $c_0({\bf L},A)$ by balls $B(c_0({\bf L},A),x,r)$ 
with $0< r\le p^{-j}$ and $x\in c_0$,
since such balls are either disjoint or one of them is contained in
another and $\Gamma _{\bf L}$ is discrete in $(0,\infty )$, where
$\Gamma _{\bf L}:=\tilde \Gamma _{\bf L}\setminus \{ 0\} $. Then
$\bigcup_{q\in J}B(c_0,x_q,r_q)$ is the clopen ball in $c_0$
with $r\le p^{-j}$, if all balls in the family $J$ have non-void
pairwise intersections. Taking $B(c_0,x,r)\cap f(X)$ we get the
statement for $f(X)$ using the transfinite sequence of the covering.
To satisfy Condition $\inf_{l\ne k}dist (B_l,B_k)>0$
we take $r_l>r_0>0$ for each ball $B_l$, where $r_0$ is a 
fixed positive number.
\par {\bf 2.7. Note.} A simplex $s$ in $\bf R^n$ may be taken 
with the help of linear functionals, for example, $\{ e_j:
j=0,...,n \} $, where $e_j=(0,...,0,1,0,...,0)$ with $1$ in 
the $j$-th place for $j>0$ and $e_0=e_1+...+e_n$, $s:=\{
x\in {\bf R^n}: e_j(x)\in [0,1] \mbox{ for }j=0,1,...,n \} $.
In the case of $\bf L^n$ or $c_0({\bf L},A)$, 
if to take $B({\bf L},0,1)$ instead of
$[0,1]$, then conditions $x_j:=e_j(x)\in B({\bf L},0,1)$ for
$j=1,...,n$ imply $e_0(x)=x_1+...+x_n\in B({\bf L},0,1)$ 
or $\sum_{j\in \lambda }e_j(x)\in B({\bf L},0,1)$ 
for each finite subset $\lambda $ in $A$ respectively
due to the ultrametric inequality (since $B({\bf L},0,1,)$ is the
additive group ), that is $s=B({\bf L^n},0,1)$
or $s=B(c_0({\bf L},A),0,1)$ correspondingly. 
Moreover,
its topological border is empty $Fr (s)=\emptyset $ and $Ind(Fr(s))=-1$.
Let us denote by $\pi _{\bf L}$ an element from $\bf L$ such that
$B({\bf L},0,1^{-}):=\{ x\in {\bf L}: |x|_{\bf L}<1\} =
\pi _{\bf L}B({\bf L},0,1)$ and $|\pi _{\bf L}|_{\bf L}=\sup_{b\in
\Gamma _{\bf L}, b<1}b=:b_{\bf L}$.
\par {\bf 2.8. Definitions. (1).} A subset $P$ in $c_0({\bf L},A)$
is called a polyhedron if it is a disjoint union of simplexes
$s_j$, $P=\bigcup_{j\in F}s_j$, where $F$ is a set, $s_j=B(c_0({\bf L},A'),
x,r)=x+\pi _{\bf L}^kB(c_0({\bf L},A'),0,1)$ are the clopen balls in $c_0(
{\bf L},A')$, $A'\subset A$, $r=b_{\bf L}^k$, $k\in \bf Z$. For each $\bf L$
we fix $\pi _{\bf L}$ and such affine transformations. The polyhedron
$P$ is called uniform if it satisfies conditions $(i,ii)$:
\par $(i)$ $\sup_{i\in F}diam(s_i)<\infty $,
\par $(ii)$ $\inf_{i\ne j}dist(s_i,s_j)>0$, where $dist (s,q):=
\inf_{x\in s, y\in q}\rho (x,y)$. 
\par By vertexes of the simplex $s=B(c_0
({\bf L},A),0,1)$
we call points $x=(x_j)\in c_0({\bf L},A)$ such that 
$x_j=0$ or $x_j=1$ for each $j\in A$, $dim_{\bf L}(s):=card(A)$.
For each $E\subset A$, $E\ne A$
and a vertex $e$ by verge $v$ of the simplex $s$
we call a subset $e+B(c_0({\bf L},E),0,1)\subset s$. 
A mapping $f: C\to D$ is called affine,
if $f(ax+by)=af(x)+bf(y)$ for each $a, b\in B({\bf L},0,1)$, 
$x, y \in C$, where $C$ and $D$ are absolutely convex subsets of 
vector spaces $X$ and $Y$ over $\bf L$.
For an arbitrary
simplex its verges and vertexes are defined with the help
of affine transformation as images of verges and vertexes of the
unit simplex $B(c_0({\bf L},A'),0,1)$.
\par Then in analogy with the classical case \cite{steil}
there are naturally defined
notions of a simplexial complex $K$ and his space $|K|$, also
a subcomplex and a simplexial mapping. The latter is characterized
by the following conditions: its rectrictions
on each simplex of polyhedra that are affine mappings
over the field $\bf L$ and images of vertexes are vertexes.
A simplexial complex $K$ is a collection of verges $v$ of a simplex 
$s$ such that each verge $v'$ of a simplex $v$ is in $K$. 
The space $|K|$ of $K$ is a subset of $|s|$ consisting of 
those points which belong to simplexes of $K$, 
where $|s|$ is a topological space of $s$ in the topology inherited from
$c_0({\bf L},A)$. This defines topology of $|K|$ also.
A simplexial complex has $dim_{\bf L}K:=\sup_{v\in K} dim_{\bf L}v$, where
$v$ are simplexes of $K$. A subcomplex $L$ of $K$ is a subcollection
of the simplexes of $K$ such that each verge of $v$ 
of a simplex $s$ in $L$ is also in $L$,
so $L$ is a simplexial complex. If $e$ is a vertex of $K$, the open star of
$e$ is the subset $St_K(e)$ of $K$ which is the union of all simplexes
of $K$ having $e$ as a vertex.
\par Instead of the barycentric subdivision in the classical case
we introduce a $p^j$-subdivision of simplexes and polyhedra for $j\in \bf N$
and for the field $\bf L$, that is a partition of each simplex
$B(c_0({\bf L},A'),x,r)$ into the disjoint union of simplexes
with diameters equal to $rp^{-j}$, where $p=b_{\bf L}^{-1}$. 
Each simplex $s$ with $dim_{\bf L}s=
card(A')$ may be considered also in $c_0({\bf L},A)$, where $A'\subset A$,
since there exists the isometrical embedding $c_0({\bf L},A')
\hookrightarrow c_0({\bf L},A)$ and the projector $\pi :
c_0({\bf L},A)\to c_0({\bf L},A')$. By a dimension of a polyhedron
$P$ we call $dim_{\bf L}P:=\sup_{(s\subset P, s\mbox{ is a simplex })}
dim_{\bf L}s$. The polyhedron $P$ is called locally finite-dimensional
if all simplexes $s\subset P$ are finite dimensional over $\bf L$,
that is, $dim_{\bf L}s\in \bf N$. For a simplex $s=B(c_0({\bf L},
A'),x,r)$ by a $\bf L$-border $\partial s$ we call the union of all
its verges $q$ with the codimension over $\bf L$ equal to $1$
in $c_0({\bf L},A')=:X$, that is, $q=e+B'$, where $B'$ are balls in
$c_0({\bf L},A")=:Y$, $(A'\setminus A")$ is a singleton, $A"\subset A'$,
$X\ominus Y=\bf L$, $Y\hookrightarrow X$. For the polyhedron
$P=\bigcup_{j\in F}s_j$ by the $\bf L$-border we call $\partial P=
\bigcup_{j\in F}\partial s_j$, where $F$ is the set.
\par {\bf (2).} A continuous mapping $f$ of a set $M$,
$M\subset c_0({\bf L},A)$ into a polyhedron $P$ we call essential,
if there is not any continuous mapping $g: M\to P$ for which are satisfied
the following conditions:
\par $(i)$ $g(M)$ does not contain $P$;
\par $(ii)$ there exists $M_0\subset M$, $M_0\ne M$ with $f(M)\cap
\partial P=f(M_0)=g(M_0)\subset \partial P$ and their restrictions
coincide $f|_{M_0}=g|_{M_0}$;
\par $(iii)$ if $f$ is linear on $[x,y]:=\{ tx+(1-t)y| t\in B({\bf L},0,1)\}
\subset M$, then $g$ is also linear on $[x,y]$ such that
$g(x)\ne g(y)$ for each $f(x)\ne f(y)$.
\par {\bf (3).} The function $f$ from \S 2.8.(2) is inessential,
if there exists such $g$.
\par {\bf (4).} Let $f: M\to N$ be a continuous mapping, $c_0({\bf L},A)
\supset N\supset P$, $P$ be a polyhedron. Then $P$ is called essentially
(or inessentially) covered by $N$ under the mapping $f$, if
$f|_{f^{-1}(P)}$ is essential (or inessential respectively).
\par {\bf (5).} Let $f: M\to P$ and $g: M\to P$ are continuous mappings,
where $M$ is a set, $P$ is a polyhedron. Then $g$ is called a permissible
modification of $f$, if three conditions are satisfied:
\par $(i)$ from $a\in M$ and $f(a)\in s$ it follows $g(a)\in s$, where
$s$ is a simplex from $P$;
\par $(ii)$ if $x$ and $y\in M$, $[x,y]\subset M$ and $f: [x,y]\to P$ is
linear, then $g: [x,y]\to P$ is also linear (over $\bf L$)
and $g(x)\ne g(y)$ for each $f(x)\ne f(y)$;
\par $(iii)$ $f(\partial M)=g(\partial M)$, when $M$ 
contains a polyhedron $Q$ with $dim_{\bf L}Q>0$ such that $Q$
is clopen in $M$.
\par {\bf (6).} The mapping $f: M\to P$ is called reducible (or irreducible),
when it may (or not respectively) have the permissible modification
$g$ such that $f(M)$ is not the subset of $g(M)$.
\par {\bf (7).} A mapping $f: P\to Q$ for polyhedra $P$ and $Q$
is called normal, if 
\par $(i)$ $\rho _Q(f(x),f(y))\le \rho _P(x,y)$ for each $x$ and $y\in P$,
that is $f$ is a non-stretching mapping;
\par $(ii)$ there exists a $p^j$-subdivision $Q'$ of polyhedra $Q$ 
such that $f: P\to Q'$ is a simplexial mapping, that is, $f|_s$ is 
affine on each simplex $s\subset P$ and $f(s)$ is a simplex from $Q'$.
\par {\bf (8).} Let $X=\lim_j \{ X_j, f^j_i, E\} $
be an expansion of $X$ into the limit of inverse spectra of polyhedra 
$X_j$ over $\bf L$. This expansion is called 
\par $(a)$ irreducible,
if for each open $V\subset X$ there exists a cofinal subset $E_V\subset E$
such that $\{ X_j, f^j_i, E_V\} $ is also the
irreducible polyhedral representation
of the space $V$, that is $f^j_i: X_j\to X_i$ are irreducible and
surjective for each $i\ge j\in E_V$. The polyhedral system (representation)
$\{ X_j, f^j_i, E \} $ is called 
\par $(b)$ $n$-dimensional, if $dim_{\bf L}X_j
\le n$ for each $j\in E$, $\sup_{j\in E}dim_{\bf L}X_j=n$,
where $n$ is a cardinal number.
The polyhedral spectrum is called 
\par $(c)$ non-degenerate, if
projections $f^j_i$ are non-degenerate, that is
$f^j_i(s)$ is not a singleton for each simplex $s$ of $X_j$.
\par {\bf 2.9. Notes.} Conditions $2.8.(2(iii), (5(ii, iii))$ and restrictions
on $f^j_i$ in $2.8.(8(a))$ are additional in comparison with the classical
case. They are imposed in \S 2.8, since there exists 
a continuous non-linear retraction 
$r: s_j\to \partial s_j$ for the non-Archimedean field $\bf L$, which may 
be constructed with the help of a $p$-subdivision
and projections of $c_0({\bf L},A)$ are on its subspaces 
$c_0({\bf L},A')$ associated with the standard basis 
$\{ e_j:$ $j\in A \} $, where $A'\subset A$.
If $f$ is simplexial on each polyhedron $M$ and $2.8.(2(iii))$ is accomplished,
then $dim_{\bf L}g(M)=dim_{\bf L}f(M)$.
\par In \cite{isb4,isb5} were studied uniform spaces and ANRU. Here we mean 
by ANRU an ultraunifrom space $X$ such that under embedding into an 
ultrauniform space $Y$ there exists 
the uniformly continuous retract $r: V\to X$
of its uniform neighbourhood $V$, $X\subset V\subset Y$. 
We denote by $U(X,Y)$ for two ultrauniform spaces $X$ and $Y$
an ultrauniform space of uniformly continuous mappings $f: X\to Y$
with the uniformity generated by a base of the form $W=\{ (f,g)|
(f(x),g(x))\in V$ $\mbox{for each}$ $x\in X\} $, where $V\in \sf V$,
$\sf V$ is a uniformity on $Y$ corresponding to ${\bf P}_Y$.
\par Here we call an ultrauniform space $Y$ injective, if for each
ultrauniform space $X$ with $H\subset X$ and an uniformly continuous mapping
$f: H\to Y$ there exists an uniformly continuous extension $f: X\to Y$.
\par {\bf 2.10. Theorem.} {\it $B(c_0({\bf L},A),0,1)$ is an injective space 
for $card(A)<\aleph _0$.}
\par {\bf Proof} may be done analogously to Theorem 9 on p. 40 
in \cite{isb5}.
\par {\bf 2.11. Theorem.} {\it Each uniform polyhedron $P$ over $\bf L$ is
ANRU.}
\par {\bf Proof.} Since $a=\inf_{i\ne j\in F}dist(s_i,s_j)>0$,
$b=\sup_{i\in F}diam(s_i)<\infty $, then there exists a 
minimal uniform covering
$\sf U$ such that for each $s_i$ there exists $a/p$-clopen neighbourhood
$V_i\in \sf U$, where $V_i=s_i^{\epsilon }$, $\epsilon =a/p$,
$A^{\epsilon }= \{ y\in Y:$ 
$d_Y(y,A) = \inf_{a\in A}d_Y(y,a)\le \epsilon \} $
for an ultrametric space $(Y,d_Y)$ is an 
$\epsilon $-enlargement of a subset $A$ in $Y$.
In this case $P$ is a subset of a Banach space $c_0({\bf L},A)$, 
where $card (A)=w(P)$.
Each $s_i$ is an uniform retract $V_i$, $r_i: V_i\to s_i$,
consequently, there exists uniformly continuous retraction $r: S=
\bigcup_{V_i\in \sf U}V_i\to P$ such that $r|_{s_i}=r_i$ for each $i$,
since $\sup_{i\in F}diam(s_i)<\infty $
and $r_i\circ r_i=r_i$, $r_i|_{s_i}=id$, 
hence $r\circ r=r$ and $r|_P=id$.
The neighbourhood $S$ is uniform, since $St(P,{\sf U}):=\bigcup_{(A\in 
{\sf U}, A\cap P\ne \emptyset )}A=S$ \cite{eng}.
\par {\bf 2.12. Note.} Further for uniform spaces are considered
uniformly continuous mappings and uniform polyhedra and for topological 
spaces continuous mappings and polyhedra if it is not specially
outlined.
\par {\bf 2.13. Lemma.} {\it Let $(X,{\bf P})$ be an ultrauniform
strongly zero-dimensional space, $P$ be a polyhedron over $\bf L$,
$A_1$,...,$A_q$ are non-intersecting closed subsets in $X$,
$q\in \bf N$, $card(P)\ge q$. Then there exists an uniformly
continuous mapping $f: X\to P$ such that $f(A_i)\cap f(A_j)=\emptyset $
for each $i\ne j$.}
\par {\bf Proof.} There exists a disjoint clopen covering $V_j$ for $X$
satisfying $A_j\subset V_j$ for each $j=1,...,q$ (see Theorem 6.2.4
\cite{eng}). Then we can take locally constant or locally
affine mapping $f$ with $f(V_j)\subset s_j\subset P$.
\par {\bf 2.14. Lemma.} {\it Each non-stretching mapping $f: E\to P$
has non-stretching continuation on $(\tilde X, \rho )$, where $P$
is a uniform polyhedron over $\bf L$, $E\subset X$,
$\tilde X$ is a completion of an ultrametric space $(X,\rho )$.}
\par {\bf Proof.} There exists an embedding $P\hookrightarrow
c_0({\bf L},A)$ for $card(A)=w(P)$ due to Lemma 4 \cite{lu3}
or Theorem 2.4 above.
We choose $f: E\to c_0({\bf L},A)$ with $f=(f^i: i\in A)$,
$f^i: E\to {\bf L}e_i$, where $(e_i: i\in A)$ is the orthonormal basis
in $c_0({\bf L},A)$ and $\inf _{(i\ne j, s_i\mbox{ and }s_j\subset P)}
dist(s_i,s_j)>0$, $s_i$ are simplexes from $P$. 
Each mapping $f^i: E\to {\bf L}e_i$ is non-stretching, 
that is $\| f^i(x)-f^i(y) \| \le \rho (x,y)$ for each $x, y \in 
E$. Evidently, $f$ has non-stretching extension on $\tilde E$ 
such that $\tilde E$ is closed $\tilde X$. On the other hand, 
due to Theorem 2.4 there exists $\epsilon >0$ such that $f$ 
has a non-stretching extension $f: E^{\epsilon }\to P^{\epsilon }$. 
By Theorem 2.11 for sufficiently small $\epsilon >0$ there exists 
a retraction $r: P^{\epsilon }\to P$ such that 
$r\circ f: E^{\epsilon }\to P$ is is non-stretching, 
$r\circ f=f$ on $\tilde E$. There 
exists $V$ clopen in $\tilde X$ such that $E\subset V\subset E^{\epsilon }$
due to Lemma 2.2, hence $r\circ f|_V$ has a non-stretching extension
of $f$ from $E$, since $(\tilde X,\rho )$ is strongly zero-dimensional.
This follows also from the fact that a free Banach space 
$B(E,\rho ',{\bf L})$ over $\bf L$ generated by $E$ has the isometric 
embedding into $B(X,\rho ',{\bf L})$, where $\rho '$ is from Lemma 2.2
(see also \cite{lu3}).
\par {\bf 2.15. Definitions.} An ultrauniform space $(X,{\bf P})$
is called $LE$-space, if each uniformly continuous $f: Y\to \bf L$ has 
an uniformly continuous extension on all $X$, where $Y\subset X$.
\par Let $X$ be a Tychonoff space. If $X$ has the finest 
uniformity compatible with its topology, then it is called fine.
\par {\bf 2.16. Theorem.} {\it An ultrametric space $X$ is an $LE$-space 
if and only if ${\tilde X}=\lim \{ X_m, f^m_n, E\} $, where $X_m$ are
fine spaces and bonding mappings
$f^m_n: X_m\to X_n$ are uniformly continuous for 
each $m\ge n\in E$, $\{ X_m, f^m_n, E \} $ is an inverse mapping system.}
\par {\bf Proof.} Let us consider $X_A=B(c_0({\bf L},A),0,1)$ for
$card(A)\ge \aleph _0$ or $X_A=\bf L^n$ which are not fine spaces, since 
on $X_A$ there are continuous $f: X_A\to \bf L$ which are not uniformly
continuous. If $X_A$ is fine, then each continuous 
function $f: X_A\to \bf L$ is uniformly continuous.
Then for the embedding of the non-fine space 
$X_A\hookrightarrow c_0({\bf L},A)$ there is not compact $H$ with 
$X_A\subset H \subset c_0({\bf L},A)$. Therefore, in $X_A$ there exists a
countable closed subset of isolated points $Y=\{ x_i: i\in {\bf N} \} $ and 
a continuous function $f: Y\to \bf L$ with $|f(x_j)-f(x_i)|/|x_j-x_i|>c_i>0$
for each $j\ge i$ and $\lim_{i\to \infty }c_i=\infty $ such that $f$ has
a continuous extension $g$ on $X_A$. Indeed, there exists a clopen 
neighbourhood $W$, that is a $b$-enlargement (with $b>0$) relative to 
the ultrametric $\rho _A$ in $X_A$, so there is a continuous retraction
$r: W\to Y$, $g(x)=const \in \bf L$ on $(X_A\setminus W)$, $g(x)=f(r(x))$
on $W$, $g|_Y=f$. Therefore, $g$ is continuous and uniformly continuous on
$X_A$ (see also Theorem 1.7 \cite{coris} ).
\section{Polyhedral expansions.}
\par {\bf 3.1. Theorem.} {\it Each complete ultrauniform space $(Y,{\bf P})$ 
is a limit of an inverse spectra of ANRU $Y_j$, where $Y_j$ are embedded
into complete locally $\bf L$-convex spaces.}
\par {\bf Proof.} In view of Corollary 2.5 there exists a uniform embedding
$Y\hookrightarrow \prod_{\rho \in \bf P}c_0({\bf L},A_{\rho })=:H$.
In each $c_0({\bf L},A_{\rho })$ may be taken the orthonormal basis
$\{ e_{j,\rho }: j\in A_{\rho } \} $, $card (A_{\rho })=w(Y_{\rho })$
and define canonical neighbourhoods $U(f,b;(j_1,\rho _1),...,(j_n,\rho _n))
:=\{ q\in H: |\pi _{j_i,\rho _i}(q)-\pi _{j_i,\rho _i}(f)|<b$ $\mbox{ for }$
$i=1,...,n \} $, where $\pi _{j,\rho }: H\to {\bf L}e_{j,\rho }$ are 
projectors, $f\in H$, $b>0$, $n\in \bf N$. Each clopen subset
$Z_b:=H\setminus U(f,b;(j_1,\rho _1),...,(j_n,\rho _n))$ is uniformly 
continuous (non-stretching) retract of $Z_{b/p}$, that is $Z_b$ is ANRU.
Analogously each finite intersections $\bigcap_{l=1}^qZ(f_l,b_l,(j_1^l,
\rho _1^l),...,(j_n^l,\rho _n^l))=\bigcap_{l=1}^kZ_{l,b_l}$ 
are also ANRU, since $Z_{k,b_k}\supset \bigcap_{l=1}^kZ_{l,b_l}$
and a retraction
$g_k: Z_{k,b_k/p} \to Z_{k,b_k}$ produces non-stretching relative to the 
corresponding pseudoultrametric retraction $g_k: \bigcap_{l=1}^k
Z_{l,b_l/p}\to \bigcap_{l=1}^kZ_{l,b_l}$. Hence $g_q\circ ... \circ g_1=g:
\bigcap_{l=1}^qZ_{l,b_l/p}\to \bigcap_{l=1}^qZ_{l,b_l}$ is the 
(non-stretching) retraction. Let $S$ be the ordered family of such
$\bigcap_{l=1}^qZ_{l,b_l}\subset H\setminus Y$, then $\bigcap_{Z\in S}
(H\setminus Z)=Y$. Further as in the proof of Theorem 7.1 \cite{isb4}.
\par {\bf 3.2. Lemma.} {\it Let $(X, \rho _X)$ and $(Y, \rho _Y)$ 
be ultrametric
spaces, $f: X\to Y$ be a continuous mapping such that $f|_H$ is a 
$b$-mapping (that is $\rho _Y(f(x),i(x))\le b$ for each 
$x\in H$), where $H$ is dense in $X$,
$b>0$, $i: H\hookrightarrow Y$ is an embedding. 
Then $X$ and $Y$ may be embedded
into a Banach space $W$ over $\bf L$ such that $X$ with equivalent 
ultrametric in $W$ and $\| f(x)-x\| \le b$ for each $x\in X$ (that is
$f$ is a $b$-mapping on $X$).}
\par {\bf Proof.} From Theorem 2.4, Lemma 2.3 and Corollary 2.6
it follows that there exists such embedding of $X$ into the corresponding
$W=c_0({\bf L},A)$ with $card(A)=w(X)$ with the disjoint clopen covering
$V=\{ B(c_0({\bf L},A),x_j,b_j)=:B_j: x_j\in H, \infty >b_j>0, j\in F \} $.
Let $Y_j:=f(B_j)\subset Y$, consequently, $diam(Y_j)\le b$, since $f$ 
is the $b$-mapping, that is $f$ realizes the covering of $X$ consisting
from subsets of diameters not greater than $b$. Then $\rho _X(x',x")\le 
\max (\rho (x',x_j), \rho (x_j, x"))$, where $x"\in f^{-1}(y)$, $y\in Y_j$,
$x'\in B(c_0({\bf L},A),x_j,b_j)$, $f(x')=y$, $diam(f^{-1}(y))\le b$.
Let $x_i\ne x_j$, this is equivalent to $B_i\cap B_j=\emptyset $ and 
is equivalent to $\rho _X(x_i,x_j)>b$, consequently, $Y_i\cap Y_j
=\emptyset $ and $Y=\bigcup_{j\in F}Y_j$. In view of continuity of
$f$ there exists the embedding of $Y_j$ into $B_j$, since $f$ is 
the $b$-mapping, $w(Y)\le w(X)$, $0<b_j\le b$ and due to Theorem 2.4.
\par {\bf 3.3.1. Lemma.} {\it Let there exists a non-stretching 
(uniformly continuous) mapping $f: R\to P$ and a non-stretching
(uniformly continuous) permissible modification $g: M\to P$, where
$R$ is a complete ultrametric space ($LE$-space respectively), 
$P$ is a polyhedron, $M$ is a subspace in $R$; if $R$ is a polyhedron let
$M$ be a subpolyhedron. Then $g$ has a non-stretching (uniformly
continuous respectively) extension on the entire $R$ and this extension is
a permissible modification of $f$.}
\par {\bf Proof.} For a complete ultrametric space $R$ the mapping
$g$ has the non-stretching (uniformly continuous) extension on
the completion of $M$ which coincides with the closure $\bar M$ of $M$
in $R$, since $R$ and $P$ are complete. The space $P$ is complete, since
it is ANRU by Theorem 2.11 (see also Theorem 1.7 \cite{isb2},
Theorems 8.3.6 and 8.3.10 \cite{eng}). For the embedding $R\hookrightarrow
c_0({\bf L},A)$ with $card (A)=w(R)$ by Theorem 2.4 it may be taken due 
to Corollary 2.6 the disjoint clopen covering $V$ such that each $W\in V$
has the form $W=R\cap B(c_0({\bf L},A),x,r_W)$, where $r_W>0$, $x\in R$,
$\sup_{W\in V}r_W\le p^{-j}$ (for each $j\in \bf N$ there exists
such $V$).
\par In view of uniform continuity of $f$ and uniformity of the polyhedron
$P$ there exists $V$ such that for each $W$ from $V$ there exists a simplex
$T\subset P$ with $f(W)\subset T$. The space $c_0({\bf L},A)$
has the orthonormal basis $\{ e_j: j\in A\} $. If $f$ is linear on no any
$[a,b]\subset R$, then for the construction of $g$ may be used Lemma 2.14
or Theorem 2.16. This may be done with the help of transfinite induction
by the cardinality of sets of vertexes of simplexes from $P$ (or using
the Teichm\"uller-Tukey Lemma), since $P$ is a disjoint union of simplexes. 
In general let $g(\bar M)$ contains
each zero-dimensional over $\bf L$
simplex $T^0\subset P$. If it is not so, then 
there exists a point $b=f^{-1}(T^0)$ in which $g$ is not defined.
If $f$ is non-linear on each $[a,b]\subset R$ with $a\in \bar M$,
then $g(b)=f(b)$. If $f$ is linear on such $[a,b]$, then 
from $\rho (a,b)<\infty $ it follows that $[a,b]$ is homeomorphic
with the compact ball $B({\bf L},a,\rho (b,a))$ in $\bf L$
and $f([a,b])$ is homeomorphic with $B({\bf L},f(a), \| f(b)-f(a)\| )$. Then
in $\bar M\cap [a,b]=Y$ there exists $y$ such that 
$\rho (y,a)=\sup_{x\in Y}\rho (a,x)=t$. The subspace $Y$ is covered 
by the finite number of $W\in V$. For $y\ne a$, that is $t>0$, $f([a,b])$
is compact and is contained in the finite number of simplexes from $P$, 
consequently, there exists the permissible modification of $g$ on 
$[a,b]$ also. 
\par Let $E$ and $F$ be two subsets in $R$. We denote by $sp(E,F,f)$
the subspace $cl((\bigcup_{(a\in E, b\in F, f|_{[a,b]}\mbox{ is }
{\bf L}-\mbox{ linear })}[a,b])\cup E\cup F)$ in $R$, where $cl(S)$ 
denotes the closure of $S$ in $R$ for $S\subset R$. If $B=\{ q:
f^{-1}(T^0)=q, T^0\subset P, g$ $\mbox{is not defined in } q \} $, then
by the Teichm\"uller-Tukey Lemma there exists the extension of $g$ on 
$sp(\bar M,B,f)=:M_0$. Let $M_j:=sp(\bigcup_{i<j}M_i,B,f)$,
where $B=\bigcup_{(T^j\mbox{ is not the subset of }g(M))}T^j$, 
$T^j$ are simplexes from $P$
with the cardianlity of sets of vertexes equal to $j$, where $j\le w(P)$.
From Lemma 3.2 it follows that conditions $2.8.(5(i-iii))$ may be satisfied
on $M_j\cap R$. Considering vertexes of $s$ from $R\cap f^{-1}(T^j)\setminus
\bigcup_{i<j}M_i$ we construct $g$ on $M_j\subset R$. Since $R\hookrightarrow
c_0({\bf L},A)$, then $\sup_{j}M_j=R$, where $M_j$ are ordered by inclusion:
$M_j\supset M_i$ for each $i\le j$.
\par {\bf 3.3.2. Note.} Henceforth, we assume that for a 
uniformly continuous mapping $f: Y\to P$
are satisfied conditions: $P$ is a uniform polyhedron,
bonding mappings $f^m_n: X_m\to X_n$ are uniformly continuous
(see \S \S 2.8, 2.16), where $Y\subset X$, $(X,\rho )$
is an ultrametric space.
\par {\bf 3.4. Lemma.} {\it Let $f: M\to P$ is irreducible, $M$ and $P$
are polyhedrons, $N$ is a subpolyhedron in $M$, $Q$ is a subpolyhedron
in $P$, $f(N)\subset Q$, then a mapping $f: N\to Q$ is irreducible.}
\par {\bf Proof.} Let $f: N\to Q$ be reducible, that is there exists
a permissible modification $g: N\to Q$ with $g(N)$ not contained in $f(N)$,
$q\in f(N)\setminus g(N)$. In view of Lemma 3.3 there exists the extension
$g: M\to P$. Let $r>0$ is sufficiently small and $U:=N_r=\{ y\in M:$
$\mbox{there exists}$ $x\in N$ $\mbox{with}$ $\rho (x,y)\le r\} $ be 
the $r$-enlargement of the subspace $N$ such that $q\notin g(N_r)$. Since $M$ 
and $P$ are (uniform) polyhedra and $f$ is (uniformly) continuous and $U$ 
is clopen in $M$, then there exists a $p^j$-subdivision $M'$ and a clopen
polyhedron $H$ in $M'$ with $N_{r/p}\subset N\subset U\subset M'$ such that
$h|_H=g|_H$, $M'\setminus H$ is the subpolyhedron in $M'$, 
$h|_{M\setminus H}=f|_{M\setminus H}$. Then $q\notin h(H)$ and 
from the irreducubility of $f$ it follows that $q\notin f(M\setminus H)$,
consequently, $h$ is the permissible modification of $f$ and $h(M)$
is not the subset
of $f(M)$. This contradiction lead to the statement of this Lemma.
\par {\bf 3.5. Lemma.} {\it Let $f: M\to P$, $M$ and $P$ be 
polyhedrons over $\bf L$. Then the condition of irreducibility 
of $f$ is equvalent to that each subpolyhedron in $Q\subset P$ is 
essentially covered.}
\par {\bf Proof.} If $f$ is irreducible, then due to Lemma 3.4 
each subpolyhedron $Q$ is essentially covered. Let vice versa  each $Q$
is essentially covered and $f$ has the permissible modification $g$ with
$P=f(M)$ not contained in $g(M)$. 
From $f(\partial M)=g(\partial M)$ and that $f$ 
is essential the contradiction with the statement $\{ P$ is 
not contained in $g(M) \} $ follows,
consequently, $f$ is irreducible.
\par {\bf 3.6. Lemma.} {\it Let $P$ and $M$ be polyhedrons.
If $f: M\to P$ has only addmissible modifications,
then $f$ is irreducible.}
\par {\bf Proof.} From $f(\partial M)=g(\partial M)$ 
and that $g$ is the permissible modification of $f$ it follows that
each subpolyhedron $Q$ from $P$ is essentially covered due to Lemma 3.3.
In view of Lemma 3.5 $f$ is irreducible.
\par {\bf 3.7. Lemma.} {\it Let there is an inverse spectra
$S=\{ R_m, f^m_n, E\} $ of polyhedra $R_m$ over $\bf L$, $f^m_n$ are
simplexial mappings, $g^l: R_l\to P$, $g^n=g^l\circ f^n_l$ for a fixed
$l$, $f_n=\lim_mf^m_n$, $g=g^l\circ f_l$, $R=\lim S$,
$f_n: R\to R_n$, $E$ is linearly ordered. If $g: R\to P$ is reducible
for a polyhedron $P$, then for almost all $n$ 
(that is, there exists $k\in E$ such that for each $n\ge k$)
$g^n: R_n\to P$ are reducible.}
\par {\bf Proof.} In view of Lemma 3.3 there exists a 
permissible modification $h$ of $g$, then $g$ and $h$ define
mappings $g^n$ and $h^n$ such that $h^n$ on $R_n$ are
permissible modifications of $g^n$. If $g$ is reducible, then by Lemma 
3.5 $g^n$ are reducible for almost all $n\in E$, since $h(R)$
is not a subset of $g(R)$ and $h(\partial R)=g(\partial R)$.
\par {\bf 3.8. Lemma.} {\it If $f: M\to N$, $g: N\to T$, $f(M)=N$,
where $M$ and $N$ are polyhedra, then from $g$ is inessential (reducible)
it follows $f\circ g$ is inessential (reducible).}
\par {\bf Proof.} If $f$ is $\bf L$-linear on $[a,b]\subset M$
and $g$ is $\bf L$-linear on $[f(a),f(b)]$, then $f\circ g$ 
is $\bf L$-linear on $[a,b]$. If $h$ is a permissible modification of
$g$, then $f\circ h$ is a permissible modification of $f\circ g$,
hence $f\circ g$ is reducible. If $g$ is inessential, then there
exists a mapping $h$ such that 
\par $(i)$ $h(N)$ does not contain $T$;
\par $(ii)$ there exists $N_0\subset N$ such that $N_0\ne N$,
$g(N)\cap \partial T=g(N_0)=h(N_0)$ and $h|_{N_0}=g|_{N_0}$;
\par $(iii)$ if $g$ is linear on $[x,y]\subset N$, 
then $h$ is linear on $[x,y]$ such that $h(x)\ne h(y)$ 
for each $g(x)\ne g(y)$. Therefore, $f\circ h: M\to T$
satisfies Conditions $(i-iii)$ for $M$ instead of $N$,
$f\circ h$ instead of $h$ and $f\circ g$ instead of $g$,
consequently, $f\circ g$ is inessential.
\par {\bf 3.9. Lemma.} {\it Suppose that there is given an irreducible
inverse mapping system $S=\{ P_m, f^m_n, E \} $ of polyhedra $P_m$ 
over $\bf L$, $M$ is closed in $P=\lim S$ such that $M_n:=f_n(M)$
are subpolyhedra in $P_n$ and for each $m>l$ permissible modifications 
$g^m_l$ for $f^m_l$ are given (on the entire $P_m$) and 
for each $n>m$ mappings 
$g^m_l\circ f^n_m$ are permissible modifications of $f^m_l$.
Then the inverse mapping system $S_M:=\{ M_n, g^m_l\circ f^n_m, E_l\}$
is irreducible, where $E_l=\{ n\in E: n\ge l \}$.}
\par {\bf Proof.} From surjectivity and
irreducibility of $f^m_l$ it follows surjectivity and irreducibility of
$g^m_l\circ f^n_m$ due to Lemmas 3.4, 3.5 and 3.8. Since $E_l$ is cofinal
with $E$, then for each $W$ open in $M$ there exists $V$ open in $P$
such that $W=M\cap V$. Let $dim_{\bf L}P_l=d$. 
We use a disjoint covering of $P_l$ by simplexes clopen in it.
There exists
$l\le m_d\in E$ such that each $d$-dimensional simplex $T^d$ from $T_l$
over $\bf L$ by means of $(f^n_l,M_n)$ is either essentially 
for each $n\ge l$ or inessentially for each $n\ge m_d$ covered.
In the second case there exists a permissible modification
$k^{m_d}_l$ such that $k^{m_d}_l[(f^{m_d}_1)^{-1}(T^d)]$ does not contain
any point from $T^d\setminus \partial T^d$, hence we get
a permissible modification $k^{m_d}_l$ of $f^{m_d}_l$
on the entire $P_{m_d}$. 
If $T^d$ is essentially covered by $(f^{m_d}_l,M_{m_d})$,
then $T^d$ is essentially covered by $(f^n,M_n)$ for each $n\ge m_d$.
\par Let $q<d$, where $q$ is a cardinal. Then there exists
$m_q\ge m_d$ such that for each simplexes $T^q$ from $P_l$ with $dim_{\bf L}
T^q=q$ is either essentially for each $n\ge l$ or inessentially
for each $n\ge m_q$ covered by $(f^n_l,M_n)$. Since $P_l$ has
a disjoint covering by such simplexes and $S$ is irreducible, then
$S_M$ is irreducible.                                   
\par {\bf 3.10. Lemma.} {\it If $T$ is a simplex from a $p^j$-subdivision
$P'$ of polyhedron $P$ then for each clopen neighbourhood
$U\supset T$ such that $U$ is a subpolyhedron in $P'$ there exists
a mapping $k: P\to P$ such that $k|_U$ is simplexial
and $k(U)=T$.}
\par {\bf Proof.} In view of Theorem 2.11 there exists the retraction
$r: P\to U$ (it is uniform if $P$ is the uniform polyhedron).
In view of Lemmas 3.2 and 3.3 there exists a simplexial mapping
$f: U\to T$ such that $f(U)=T$, hence $k=f\circ r$ is the desired mapping, 
since $r|_U=id$ and $r\circ r=r$, where $id(x)=x$ for each $x\in U$.
To construct $f$ we consider simplexes $s$ of $P'$ forming its disjoint 
clopen covering. For each such $s$ either $s\cap T=\emptyset $
or $s=T$. If $s=B(c_0({\bf L},A),x,r)$ and $T=B(c_0({\bf L},A'),y,r')$
with $A'\subset A$, then $f|_s$ is evidently defined such that
$f(s)=T$. For $A'\supset A$ and $A'\ne A$ we can map $s$ on the 
corresponding verge (or vertex) of $T$. For $s=T$ we can take $f|_T=id$.
\par {\bf 3.11. Lemma.} {\it  Let $P$ be a uniform polyhedron, $f: M\to P$ 
and $g: M\to P$ are two $b$-close mappings (that is, $\rho (f(a),g(a))\le b$
for each $a\in M$), where $b<\inf_{s\subset P}diam(s)$, $s$ are clopen 
simplexes in $P$, then there exists $k$ from Lemma 3.10 such that
$k\circ g$ and $f$ are $b$-close and $k\circ g$ is a permissible modification of
$f$, where $g|_{[x,y]}$ is $\bf L$-linear or $g(x)\ne g(y)$
if and only if $f|_{[x,y]}$ is $\bf L$-linear or $f(x)\ne f(y)$
respectively for each $[x,y]\subset M$.}
\par {\bf Proof.} Let $P'$ be a subdivision of $P$ constructed with 
the help of $p^j$-subdivisions of $P$ such that 
$j\in \bf N$.
In view of $b<\inf_{s\subset P}diam(s)$ we have that 
$g(a)\in \tau $ if and only if $f(a)\in \tau $ for each simplex 
$\tau $ in $P$ and $a\in M$.
If $M_1$ is a clopen polyhedron in $M$
and $s$ is a simplex (or its verge) in $M_1$, 
then $f|_s$ is affine if and only if
$g|_s$ is the affine
mapping, moreover $f(s)$ and $g(s)$
have the same dimension over $\bf L$ and $dim_{\bf L}g(q)=
dim_{\bf L}f(q)$ for each verge $q$ in $s$.
Therefore, for a polyhedron $M_1$ we can choose $k: P\to P$
such that $k\circ g(\partial M)=f(\partial M)$, since
$P$ is a uniform polyhedron and there exists $j\in \bf N$ such that
$p^{-j}\sup_{s\subset P}diam(s)<b$.
\par {\bf 3.12. Lemma.} {\it Let $f: P\to Q$ be a non-stretching mapping
of a uniform polyhedron $P$ onto a uniform polyhedron $Q$ over $\bf L$.
Then there exists a $p^j$-subdivision $P'$ of $P$ and a normal mapping
$g: P\to Q$ such that $g$ is a permissible modification of $f$.}
\par {\bf Proof.} For each $b>0$ 
there exist a $p^j$-subdivision $P'$ of $P$ and a $p^i$-subdivision
of $Q$ and a simplexial mapping $h: P'\to Q$, which is $b$-close to $f$,
since $f$ is non-stretching.
Indeed, simplexes $T^l\subset P'$ are disjoint and clopen for them
due to Lemmas 3 and 4 \cite{lu3} $h$ can be chosen such that:
\par $(i)$ $|h(e)-f(e)|<b$ 
for linearly independent vertexes $e_{l,j}\in T^l$,
that is, $\{ (e_{l,j}-e_{i,0}): j\in A_l \}$ are linearly independent,
$e_{l,0}$ is a marked vertex from $T^l$, $card(A_l)=dim_{\bf L}T^l$ and 
\par $(ii)$ $h(T^l)$ are simplexes in $Q$ with $diam(h(T^l))\le diam(T^l)$
for suitable $b>0$ and the $p^i$-subdivision $Q'$ of $Q$, where $h|_{T^l}$ 
are affine mappings for each $l$. This is possible due to uniformity
of polyhedra $P$ and $Q$. Taking $g=k\circ h$, where $k$ is the suitable 
mapping from Lemma 3.11 we get the desired $g$.
\par {\bf 3.13. Lemma.} {\it Let $g$ be a permissible modification of
$f: R\to P$ and $h: P\to Q$ be a normal mapping, where $P$ and $Q$ are 
uniform polyhedra. Then $h\circ g$ is the permissible modification 
of $h\circ f$.}
\par {\bf Proof.} If $f(\partial R)=g(\partial R)$, then
$h\circ f(\partial R)=h\circ g(\partial R)$. 
If $f|_{[x,y]}$ is $\bf L$-linear,
then $g|_{[x,y]}$ is $\bf L$-linear, where $[x,y]\subset R$,
hence if $h\circ f|_{[u,v]}$ is $\bf L$-linear, then $h\circ g|_{[u,v]}$
is $\bf L$-linear and $h\circ g(u)\ne h\circ g(v)$ for each
$h\circ f(u)\ne h\circ f(v)$, since $h$ is affine for some 
$p^j$-subdivision $Q'$ of $Q$. If $a\in R$ and $f(a)\in s$,
where $s$ is a simplex from $P$, then $g(a)\in s$, hence 
if $h\circ f(a)\in \tau $, then $h\circ g(a)\in \tau $,
where $\tau $ is a simplex in $Q$, since $\rho _Q(h\circ f(a),
h\circ g(a))\le \rho _P(f(a),g(a))$ and $h: P\to Q'$ is 
simplexial.
\par {\bf 3.14. Lemma.} {\it Let $\{ P_n, f^n_m, E \}=S$ be an irreducible
inverse system of uniform polyhedra over $\bf L$, $M$ be closed in
$P=\lim S$, $M_n=f_n(M)$ be subpolyhedra in $P_n$, 
$\{ Q_k, g^k_l, F \}$  be an inverse spectra of polyhedra $Q_k$,
which appear by $p^{j(k)}$-subdivisions of $P_{n(k)}$, $g^k_l$ be
normal and permissible modifications of $f^{n(k)}_{n(l)}$, where
$F$ is cofinal with $E$, $N_k=g_k(M)$, $g^k_l|_{N_k}$ and
$g_l|_M$ are irreducible. Then $N_k$ are polyhedra.}
\par {\bf Proof.} This follows from Lemmas 3.7-3.12, since 
an image of a polyhedron under simplexial mapping is a polyhedron
(see also the proof of Lemma IV.31.XII \cite{freu} and 
\cite{isb3} starting the consideration from some fixed $n_1=q\in E$).
\par {\bf 3.15. Lemma.} {\it Suppose that for a complete ultrametric space
$R$ there is a non-stretching mapping $f: R\to P$, $f(R)=P$, $P$ is a
uniform polyhedron over $\bf L$. Then for each $b>0$ there exists
a $b$-mapping $g: R\to Q$ and $Q$ is a uniform polyhedron over $\bf L$ 
such that for sufficiently fine $p^j$-subdivision $Q'$ of a polyhedron $Q$
there exists a normal mapping $k: Q\to P$ and $k\circ g$ is a non-stretching
permissible modification of $f$.}
\par {\bf Proof.} For $R$ there exists the embedding $R\hookrightarrow
c_0({\bf L},A)$ with $card (A)=w(R)$ by Theorem 2.4 and a clopen
neighbourhood $S$ with $R\subset S\subset R(r)$ that is a uniform polyhedron
due to Corollary 2.6, where $R(r)$ denotes the $r=b/p$-enlargement of 
$R$. By Lemma 3.2 there exists the 
subpolyhedron $Q$ with $R\subset Q\subset S$ 
and the $b'$-mapping $g: R\to Q$. If $[a,z]\subset S$ and $f|_{[a,z]}$ is 
$\bf L$-linear, then we can choose $g$ such that it is linear on $[a,z]$
and $g(a)\ne g(z)$ when $f(a)\ne f(z)$. From the completeness of $R$ and Lemma
2.14 it follows that $f$ has the non-stretching extension $f: S\to P$.
Then there are a $b$-mapping $g$, a uniform polyhedron $Q$ 
and a non-stretching mapping $h: Q\to P$ for sufficiently small
$r>0$ and $b'>0$. For $h$ due to Lemma 3.12 there exists the permissible 
modification and normal $k: Q'\to P$ such that $k\circ g(R)=P$, 
$f(\partial R)=k\circ g(\partial R)$ for $\partial R\ne 0$, that is, when
$R$ contains a clopen polyhedron $R_1$ with
simplexes $T$ having $dim_{\bf L}T>0$. Such $k$, $h$ and $g$
may be constructed on each simplex $T$ from $Q'$ and then on the entire
space.
\par {\bf 3.16. Lemma.} {\it Suppose that $R$ is a complete ultrametric space,
$f_n: R\to P_n$ are non-stretching $b_n$-mappings, $P_n$ are uniform polyhedra 
over $\bf L$, $n\in E$, $E$ is an ordered set such that for each
$b>0$ there exists $l\in E$ with $0<b_n<b$ for $n>l$. Then there
exists a normal irreducible inverse mapping system $S=\{ Q_m, g^n_m, F\} $,
$F$ is cofinal with $E$, $\lim S=R$, $Q_m$ are subpolyhedra
of $p^{j(m)}$-subdivisions of $P_{n(m)}$ and $g_m$ 
are permissible modifications
of $f_{n(m)}$, where $g_m=\lim_n g^n_m: R\to Q_m$.}
\par {\bf Proof.} In view of Lemmas 3.5, 3.6, 3.13 and 3.14 it is sufficient 
at first to construct $S$ with non-stretching normal normal and surjective
mappings $g^m_n$. This can be done due Lemma 3.15 with $g^m_l\circ f_{n(m)}$
being non-stretching permissible modifications of $f_{n(l)}$ and 
$\lim_m g^m_l\circ f_{n(m)}$ are non-stretching permissible modifications
of $f_{n(l)}$.
\par {\bf 3.17. Lemma.} {\it If the ultrametric space $(X,\rho )$
is isomorphic with 
$$\lim \{ (X_m, \rho _m); f^m_n; F\}$$ 
and the following conditions are satisfied:
\par $(1)$ for each $m$ there are embeddings $q_m: X_m\hookrightarrow
(E, \rho )$ into a complete space $(E, \rho )$; 
\par $(2)$ $f_m: X\to X_m$ are projections;
\par $(3)$ $(X_m, \rho _m)$ are ultrametric spaces;
\par $(4)$ there is given a family $\{ b_m>0: m\in F\}$, $b_m\in
\rho (X,X)$ and for each $b>0$ the set $\{ m: b_m>b \}$ is finite,
where $t_m:= \inf_{\rho (x,y)>b_m} \rho (q_m(x), q_m(y))$
and $\lim_m t_m=0$,
for each $m>n$ and $x$ and the inequality $\rho (q_m(x),
q_n\circ f^m_n(x))<t_n$ is accomplished.
Then the mappings $q_m\circ f_m$ converge uniformly to
the embedding $X\hookrightarrow E$.}
\par {\bf Proof.} In view of Lemma 2.2 we may suppose that $\rho (X,X)$
and $\rho _m(X_m,X_m)\subset \tilde \Gamma ({\bf L})$. If $x$ and $y\in X$
and $\rho (x,y)>b_n$, then $\rho (q_n\circ f_n(x), q_n\circ f_n(y))\ge t_n$.
From conditions $(1-4)$ it follows the existence of $k$ such that
$\rho (q_m\circ f_m(x), q_m\circ f_m(y))\ge t_n$ for each $m>k$ and 
for the ultrametric
$\rho $, consequently, $q=\lim_m q_m$ is the embedding, since 
$\rho (q(x), q(y))\ge t_n$.
\par {\bf 3.18. Theorem.} {\it Let $(X,{\bf P})$ be a complete ultrauniform 
space and $\bf L$ be a locally compact field.
Then there exists an irreducible normal expansion of $(X, {\bf P})$
into the limit of the inverse mapping system 
$$S = \{ P_n, f^n_m, E\} $$ 
of uniform polyhedra
$P_n$ over $\bf L$, moreover, $\lim S$ is isomorphic with 
$(X, {\bf P})$, in particular for an 
ultrametric space $(X, \rho )$ the system
$S$ is the inverse sequence.}
\par {\bf Proof.} From Corollary 2.5 and Theorem 3.1 it follows the existence
of the expansion of $(X, {\bf P})$ into the uniformly isomorphic limit of
the inverse mapping system $R=\{ Y_j, f^j_i, F \}$ of ANRU $Y_j$ with 
non-stretching $f^j_i$, where $Y_j$ are the complete ultrametric spaces.
Each $Y_j$ is closed in the finite products of the spaces $c_0({\bf L}, A_k)$.
From Lemmas 3.2-3.17 it follows the existence of the irreducible normal uniform
polyhedral expansion for each $Y_j$, moreover, using the permissible
modifications of $g^j_i$ of $f^j_i$ and the same Lemmas we can construct 
the system of the entire space $(X, {\bf P})$ with the same properties.
\par We consider further uniform coverings $V$ corresponding to the
uniform polyhedra $P=\bigcup_{W\in V}W$, which due to Theorem 2.11 are ANRU.
Let $\rho $ be the pseudoultrametric in $X$ and $\rho (X,X)\subset 
\tilde \Gamma ({\bf L})$. 
At first,
if $V$ is given with the help of the chosen pseudoultrametric
$\rho $, then we can associate with $V$ a $p^k$-nerve 
with $k\in \bf Z$, that is, an abstract simplexial complex $N_k$ vertexes
of which are elements from $V$. Its simplexes are the spans of (pulled on)
vertexes $W_j$ satisfying $\rho (W_j,W_i)\le p^kb$, where $b=\sup_{W\in V}
diam (W)<\infty $, each rib $[W_j, W_i]$ from $s$ has the lenghts not less
than $t=\inf_{W_j\ne W_i\in V}\rho (W_j,W_i)>0$. 
Then from $N_k\hookrightarrow
c_0({\bf L},A_k)$ with $card (A_k)=w(N_k)$ it follows that each $s$ is
uniformly isomorphic with some ball $B(c_0({\bf L},A),0,1)$, 
where $card (A)=m\le w(N_k)$. With each $V$ is associated the equivalence 
relation: $xRy$ if and only if there exists $W\in V$ such that $x$ 
and $y\in W$, then the quotient mapping 
$f: X\to X/R$ is defined (see Propositions 2.4.3 and 2.4.9 \cite{eng}). 
In particular with each locally finite functionally open covering
$V$ is associated the partition of 
the unity $\{ f_W: W\in V\}$, $f_W: X\to B({\bf L},0,1)$, $f_W(x)=1$ 
for $x\in W$ and $f_W(x)=0$ for $x\notin W$, $\{ f_W: W\in V\}$ is
subordinated to $V$ (see aslo \S 5.1 \cite{eng}). 
There are the canonical non-stretching mappings
$F_k: X\to N_k$. If $X$ is compact, then $X/R$ is the finite discrete
space and $dim_{\bf L}N_k=n\in \bf N$. Into each $V$ may be refined
a disjoint clopen uniform covering $K$ with $\sup_{W\in K}diam (W)
\le bp^{-j}$, where $j\in \bf N$. That is, $V$ has the uniform strict
shrinking. 
\par In general there exist abstract simplexial complexes
$N_V$ with $dim_{\bf L}N_V\ge dim (X)\ge 0$,
if $V$ is an arbitrary functionally open covering of $X$ of order
not less than $dim (X)$. 
For $dim(X)=\infty $ we interpret orders of coverings
as cardinals equal to $k\ge \aleph _0$.
This is due to the following second procedure.
If $\{ W_j:$ $j \in \Lambda _W \} $
are elements of a subfamily $\sf W$ of a covering $V$ such that
$W_j$ are pairwise distinct and $\bigcap_{j\in \Lambda _W}W_j\ne \emptyset $,
then to $\sf W$ corresponds an abstract simplex $s$ in $c_0({\bf L},
\Lambda _W)$ such that $W_j$ correspond to vertexes of $s$.
Evidently, if the simplex $s$ is clopen in 
the Banach space $c_0({\bf L},A)$,
then it is characterized by two points $x, y \in s$ such that
$s=B(c_0({\bf L},A),x,|x-y|)$. If $s$ with $Int(s)=\emptyset $,
then vertexes of $s$ characterize a minimal Banach subspace
of $c_0({\bf L},A)$ in which $s$ is contained.
\par Using the first procedure we 
can consider the sequence of such shrinkings: $V^{m+1}\subset
V^m$ with $b_m=bp^{-m}$, where $m\in \bf N$. With each $V^m$ is associated
$p^{k(m)}$-nerve $N_{k(m)}$. Let $k(m)\ge -m$, $k(m+1)\le k(m)$ for
each $m\in \bf N$ and $\lim_{m\to \infty }k(m)=-\infty $. If $x$ is 
an isolated point in $X$, then there exists $n\in \bf N$ 
with $\max (b_n, p^{k(n)}b_n)
<\inf_{y\in X\setminus \{ x\} } \rho (x,y)$. Then the simplex $s\subset
N_{k(m)}$ with $x\in s$ and $m\ge n$ is zero-dimensional over $\bf L$, 
that is, $s=\{ x\} $.
\par By the construction of $N_k$ for each simplex 
$s_{m+1}\subset N_{k(m+1)}$ there exists $s_m\subset N_{k(m)}$ with
$f^{m+1}_m(s_{m+1})\subset s_m$, where $f^j_i$ are the bonding
mappings of the inverse
sequence $S=\{ N_{k(m)}, f^m_i, {\bf N})$. Each $f^{m+1}_m$ is non-stretching,
since decreases the distance at least into $p$ times and $b_m/b_{m+1}\ge p$.
If $x\ne y$, then there exists $n$ with $\max (b_n, b_np^{k(n)})< \rho (x,y)$,
consequently, for each $m>n$ there are disjoint simplexes $s$ and $s'\subset
N_{k(m)}$ with $x\in s$ and $y\in s'$. Therefore, there exists the 
uniformly continuous mapping $g: X\to \lim S$, where $g(x)=\lim_m
\{ s_m, f^m_i\}$ and $s_m\ni x$ for each $m\in \bf N$. Therefore, 
the uniformly continuous projectors $f_m: X\to N_{k(m)}$ are defined,
since for each $b>0$ there exists $r\in \bf N$ 
such that $b_mp^{k(m+r)-k(m)}<b$
and $f_m(W)=f^{m+r}_m\circ f_{m+r}(W)$ and $diam (f_m(W))<b$, where
$W\in V^{m+r}$, $f_{m+r}(W)$ belongs to clopen star of the corresponding
vertex $v\in N_{k(m+r)}$. Further (see also Lemma IV.33 in \cite{isb5}) we
can verify that $g(X)=\lim S$ and $g$ is the uniform isomorphism,
since $\bigcap_{m=l}^{\infty }g_m^{-1}(s_m)= \{ x \} $
for the family $\{ s_m: m \} $ corresponding to $x$ (see above).
\par Evidently, $dim_{\bf L}N_{k(m)}$ may be from $0$ for $k(m)=-m$ up to
$card(A)$ with $card (A)=w(X)$. For $k(m)>-m$ in the inverse sequence $S$ the
mappings $f^{m+1}_m$ may map simplexes $s$ from $N_{k(m+1)}$ into simplexes
$q$ from $N_{k(m)}$ of lower dimension over $\bf L$, for example,
when $W_{m+1}\subset W_m$, $W_m\in V^m$, $W_m=B(c_0({\bf L}, A_j),x,r)$,
$W_{m+1}=B(c_0({\bf L},A_n),x',r/p)$, $card (A_n)>card(X/R_m)\ge
dim_{\bf L} N_{k(m)}\ge card (A_j)$, since $dim _{\bf L}N_{k(m+1)}\ge 
card(A_n)$ for $k(m+1)>-m-1$. This fact is based on the 
homeomorphism of $D^{\aleph _0}$, $B({\bf L},0,1)$
and $B({\bf L},0,1)^{\aleph _0}$, where $D= \{ 0,1 \} $
is a discrete two-point set (see also \cite{eng}).
\par For the complete ultrauniform space $(X, {\bf P})$ we can consider
the base of uniform coverings $\{ V^n_{\rho }: n\in {\bf N}, \rho \in
{\bf P} \} $, where each $V^n_{\rho }$ is given relative to the considered
$\rho $, $S=\{ N_{\rho , k(m)}, f^{\rho , k(m)}_{\rho ', k(m')},
{\bf P\times N} \}$. To each $V^m_{\rho }$ there 
corresponds $N_{\rho , k(m)}$;
$\rho '\le \rho $ if and only if $\rho '(x,y)\le \rho (x,y)$ for each
$x$ and $y\in X$; $(\rho ', m') \le (\rho ,m)$ if and only if
$\rho '\le \rho $ and $m'\le m$. If to associate $N_V$ with
$dim (X)$ and orders of arbitrary functionally open coverings 
$\sf {\tilde V}^n$ of the space $X$, then there exist an inverse  
spectrum of dimension equal to $dim (X)$, where
${\sf \mbox{ }^*{\tilde V}^m}< \sf {\tilde V}^n$
for each $m<n\in E$, that is $\sf {\tilde V}^m$
is a star refinement of $\sf {\tilde V}^n$.
We take $\{ {\sf {\tilde V}^m}: m\in E \} $ such that
for each $x\in X$ there exists a cofinal subset $E(x)\subset E$
and a family $ \{ {\tilde V}^m_{\beta }: \beta \in E(x) \} $
with $\bigcap_{\beta \in E(x)}{\tilde V}^m_{\beta }= \{ x \} $.
Let $P_m=N_{\sf {\tilde V}_m}$ for each $m\in E$ and
$f^m_n: P_m\to P_n$ be uniformly continuous bonding mappings,
$g(x)=\lim_m \{ s_m, f^m_i \} $,
where $s_m$ are such that $s_m\ni x$ for each $m\in E$.
Polyhedra $P_m$ can be chosen such that $\sup_{(s\subset P_m)}diam(s)<
p^{k(m)}b$ and $\lim_{m\in F}k(m)=-\infty $ for each linearly ordered 
subset $F$ of $E$ such that $F$ is cofinal with the first countable 
ordinal $\omega _0$, $s$ are simplexes in $P_m$, $k(m)\le k(n)$ for each 
$m\ge n\in F$, where $k(m)\in \bf Z$, $b=const >0$.
Then $g: X\to \lim S$ is a uniform isomorphism.
Using permissible modifications and Lemmas 3.2-3.17 we get the statement 
of this theorem.
\section{Absolute polyhedral expansions and their applications.}
\par Many of the results from \cite{koz} about absolute polyhedral 
expansions in the classical case (that is, for polyhedra over $\bf R$)
may be transformed in the non-Archimedean case. Further the main differences
are given. Part of lemmas and definitions from \cite{koz} may be naturally
reformulated and proved, for example, with the substitution of barycentric
subdivision into $p^j$-subdivision.  We use Theorems 3.1, 3.18
and apply absolute polyhedral expansions at first to
each ultrametric space $Y_j$, then to the entire ultrauniform 
space $(X,{\bf P})$. The polyhedra $P$ considered in \S 4 as well as 
in \S \S 2, 3 may have in general infinite dimension over $\bf L$.

\par {\bf 4.1.1. Definitions and Notes.} A proper $p^j$-subdivision 
of a complex $K^a$ over $\bf L$ is called normal. 
Certainly $dim_{\bf L}K^a$ may be infinite.
With the help of affine mappings
for each simplexial complex $K^a$ the space $|K^*|$
of its $p^j$-subdivision is the same as $|K^a|.$ This space is 
denoted by $K$ and is supplied with the ultrametric induced
from $|K^a|$. A support of a subset $Y\subset K$ relative to complex $K^a$
is the least subcomplex containing $Y$. A support of a point
is a simplex $s$ in $K^a$ of minimal $dim_{\bf L}s$ containing it. 
If $K^a_1$ is a subcomplex (clopen) of a complex $K^a$, then the set
$K\setminus K_1$ is called an open (clopen) subcomplex of a complex $K^a$.
If $A$ is an open (clopen) subcomplex of a complex $K^a$, then
a set $K\setminus A$ is a space of some subcomplex (clopen) 
of $K^a$, which is denoted by $K^a\setminus A$. The closure of 
an open subcomplex $A$ of a complex $K^a$ is a subspace
of all subcomplexes of $K^a$ consisting of all simplexes
interiors of which are contained in $A$ and in addition of all its 
verges in the non-Archimedean sense (see \S 2).
The least clopen subcomplex of the 
complex $K^a$ containing $X$, $X\subset K^a$, is called the polyhedral 
neighbourhood or the clopen star of a set $X$ 
relative to $K^a$ and is denoted
$K^aX$. They form a pointwise finite covering of $K$.
Moreover, a point $x\in K$ is contained in a clopen simplex
$\tau $ if and only if its vertexes $\{ e_j: j\in J \} $
are all vertexes $\{ e_j: j\in J \} $ stars of which
contain this point $x$. A clopen star of a subcomplex 
is by definition a clopen star of its space.
A clopen star $K^aK_1$ of a subcomplex $K^a_1$ in 
a complex $K^a$ is by definition a union of all 
clopen simplexes in $K^a$ which have nonvoid verge in 
the subcomplex $K^a_1$. This means that $K^aK_1$ is 
also a union of clopen stars of vertexes of 
the subcomplex $K^a_1$. 
\par Subcomplex of a complex $K^a$ a space of which is a clopen star
of a set $X$ (subcomplex $X^a$) is called a star
of the set $X$ (subcomplex $X^a$) and is denoted by
$(K^aX)^a$ or by ${\bar K}^aX^a$,
if $X^a$ is the subcomplex of the complex $K^a$.
\par {\bf 4.1.2. Lemma.} {\it Let $K^*$ be a normal, that is
multiple $p^j$-subdivision, of a complex $K^a$ and $t$ be a simplex
in $K^*$. Then there exists a vertex $e\in K^a$ such that
$K^at\subset K^ae$. Moreover, if $K^a_1$ is the subcomplex of the complex 
$K^a$ such that $|t|\cap K_1\ne \emptyset $, then $e$ 
is the vertex in the subcomplex $K^a_1$.}
\par {\bf Proof.} By the definition for each vertex $e\in t$
we have $K^ae\subset K^a_1$, if $|t|\cap K^a_1\ne \emptyset $.
Since simplexes $s$ in $K^*$ form a disjoint covering
of $|K^*|$, then $K^*t\subset K^*e$.
\par {\bf 4.1.3. Definition.} A subcomplex $K^a_1$ of a 
complex $K^a$ is called complete, if $s\subset K^a_1$
for each simpex $s$ all vertexes of which are in $K^a_1$.
\par {\bf 4.1.4. Lemmas.} {\it Let $K^a_1$ be a complete
subcomplex of a complex $K^a$.}
\par {\bf (i).} {\it If $K^*$ is a $p^j$-subdivision 
of a complex $K^a$, then a $p^j$-subdivision 
$K^*_1$ of the subcomplex $K^a_1$ induced by $K^*$ 
is the complete subcomplex of the complex $K^*$.}
\par {\bf (ii).} {\it If $f: N^a\to K^a$ is a simplexial mapping,
then $f^{-1}K_1^a$ is a complete subcomplex of the complex $N^a$.}
\par {\bf Proof.} $(i)$. If $\tau $ is a simplex in $K^a$, then
$\tau \cap K^a_1$ is a verge (may be void) of the simplex $\tau $.
On the other hand, there exists embedding of $K^a$
into $c_0({\bf L},w(|K^a|))$. Therefore, $|K^a|$
has a disjoint paving by clopen balls $B_l$ in 
$c_0({\bf L},w(|K^a|))$. For each simplex $s$ in $K^*$ there exists
a ball $B$ in $c_0({\bf L},w(|K^a|))$ such that
$B\cap |K^*|=s$. Each $p^j$-subdivision of the simplex
$\tau $ in $K^a$
consists of a disjoint union of simplexes $s_i$ in $K^*$.
\par $(ii)$. If $\tau $ is a simplex in $N^a$, then $f|_{\tau }$
is an affine mapping, hence $f(\tau )$ is a disjoint union of simplexes
from $K^a_1$. Therefore, $f^{-1}(f(\tau ))$ is covered by disjoint 
family of simplexes from $N^a$. The same is true for verges of simplexes 
from $N^a$.
\par {\bf 4.1.5. Note.} In view of Lemmas 4.1.4 above
Propositions I.1.3 and I.1.5 from \cite{koz} are true also 
in the non-Archimedean case. 
\par {\bf 4.1.6. Definition.} Subcomplexes $K_1^a$ and $K_2^a$ of
a complex $K^a$ are called adjacent, if 
$K_1=K\setminus K^aK_2$ and $K_2=K\setminus K^aK_1$.
\par {\bf 4.1.7. Lemma.} {\it Adjacent subcomplexes $K_1^a$ and
$K_2^a$ of a complex $K^a$ are complete subcomplexes.}
\par {\bf Proof.} If $\tau $ is a simplex or a verge of a simplex
in $K^a$, then either $\tau \subset K_1^a$ or $\tau \subset
K_2^a$ and $K_1^a\cap K_2^a=\emptyset $.
\par {\bf 4.1.8. Definition.} For complexes $K^a_1$ and 
$K^a_2$ their join $K^a_1\circ K^a_2$
is defined in the corresponding Banach subspace $Z$ of
$c_0({\bf L},A_1)\oplus c_0({\bf L},A_2)\oplus {\bf L}=
c_0({\bf L},A_1+A_2+1)$
with $card (A_j)=w(K^a_j)$ analogously to the classical case with substitution
of $\{ x: 0\le x \le 1 \} =[0,1]_{\bf R}\subset \bf R$ into $B({\bf L},0,1)
=[0,1]\subset \bf L$, where $j=1$ or $j=2$. 
Each simplexes $\tau _1$ in $K^a_1$ and $\tau _2$ 
in $K^a_2$ are characterized completely by their
vertexes. Then each simplex $\tau $ in $K^a_1\circ K^a_2$
has a set of vertexes $\Lambda _{\tau }$, which is a union of 
sets $\Lambda _{\tau _j}$ of vertexes of the corresponding simplexes
$\tau _j$ in $K^a_j$, where $j=1$ or $j=2$.
When $K^a_1\cap K^a_2=\emptyset $, then $Z=c_0({\bf L},A_1+A_2+1)$.
In general $Z=cl[sp_{\bf L}(K_1\cup K_2)]$.
\par {\bf 4.1.9. Lemma.} {\it $K^a_1$ and $K^a_2$ are adjacent
subcomplexes of a complex $K^a$ if and only
if each simplex $\tau $ of a complex $K^a$ has the form
\par $(i)$ $\tau =\tau _1\circ \tau _2$, where $\tau _j$ are simplexes in
$K^a_j$, $\tau _1\cap \tau _2=\emptyset $, $j=1$ or $j=2$.}
\par {\bf Proof.} If subcomplexes are adjacent, then
$\tau _1\circ \tau _2$ is obtained by junction of each
pair of points $v_1\in \tau _1$ and $v_2\in \tau _2$ by 
a segment $[v_1,v_2]$ in the sense of \S 2.8.(2).(iii).
On the other hand, equation $(i)$ is equivalent to
$\Lambda _{\tau _1}\cap \Lambda _{\tau _2}=\emptyset $.
The last equation for each pair of $\tau _1$ in $K^a_1$
and $\tau _2$ in $K^a_2$ implies that these subcomplexes 
are adjacent.
\par {\bf 4.1.10. Note.} In view of Lemma 4.1.9 above
Propositions I.1.10 and I.1.11 from \cite{koz} are true also in the 
non-Archimedean case.
\par {\bf 4.1.11. Definition and Note.} Let $K^a_1$ and $K^a_2$
be adjacent subcomplexes of a complex $K^a$ and $K^*_1$
be a $p^j$-subdivision of a complex $K^a_1$.
For each simplex $\tau _1$ of a complex $K^a_1$ let $\tau ^*_1$ 
denotes be its $p^j$-subdivision induced by $K^*_1$.
For each simplex $\tau $ in $K^a$ let $\tau ^*=\tau ^*_1\circ \tau _2$,
where $\tau _1=\tau \cap K^a_1$, $\tau _2=\tau \cap K^a_2$.
Then the union $Sd_{K^*_1}K^a$ of complexes 
$\tau ^*$ for all simplexes
$\tau $ of $K^a$ is a subdivision of a complex $K^a$.
This gives possibility to consider more general subdivisions
than $p^j$-subdivisions.
\par Evidently, analogs of Propositions I.1.12, I.1.13 and I.1.15
from \cite{koz} are true in the non-Archimedean case.
The subdivision $Sd_{K^*_1}K^a$ of the complex $K^a$ is called 
the minimal continuation of a $p^j$-subdivision $K^*_1$ of a complete
subcomplex $K^a_1$ of the complex $K^a$.
\par {\bf 4.1.12. Lemma.} {\it $Sd_{K^*_1}K^a=K^a_2\cup Sd_{K^*_1}
(K^aK^a_1)$, where $K^aK^a_1$ is a clopen star of 
the subcomplex $K^a_1$.}
\par {\bf Proof.} $K\setminus (K_1\cup K_2)$
is an intersection of clopen stars of complexes $K^a_1$ and $K^a_2$.
\par {\bf 4.1.13. Definitions and Notes.}
A mapping $h: K\to N$ is called linear relative to 
complexes $K^a$ and $N^a$, if it is simplexial for the corresponding 
polyhedra. Since $K^aS$ is the least clopen subcomplex of a 
complex $K^a$ containing a subset $S\subset K$, then $K^aS\subset 
h^{-1}(N^ah(S))$, that is equivalent to $h(K^aS)\subset N^ah(S)$.
For each mapping $f: X\to K^a$ by $\omega _fK^a$ is denoted a
covering of a space $X$ by $f^{-1}K^ae$, where $K^ae$ are called 
main stars for vertexes $e$ in $K^a$. 
Let $<N^a,K^a>_f$ denotes a subcomplex of a join $N^a\circ K^a$
of the complexes $N^a$ and $K^a$ such that its simplexes
are simplexes of $N^a$ and $K^a$ and also
$f(\tau )\circ \tau $ and all its verges, where $\tau $ is a simplex
of the complex $K^a$. Then $<N^a,K^a>_f$ is called a cylinder
of the simplexial mapping $f$. There exists a retraction
${\tilde f}: <N^a,K^a>_f \to N^a$ such that
${\tilde f}(x)=x$ for each $x\in N$, ${\tilde f}(x)=f(x)$ for 
each $x\in K$.
Evidently, Propositions I.2.1, I.2.4-10
from \cite{koz} are true in the non-Archimedean case also.
\par {\bf 4.1.14. Lemma.}  {\it $dim_{\bf L}<N^a,K^a>_f=
\max (dim_{\bf L}N^a, dim_{\bf L}f(K^a)+dim_{\bf L}K^a+1)$.}
\par {\bf Proof.} $dim_{\bf L}f(\tau )\circ \tau \le dim_{\bf L}f(\tau )
+dim_{\bf L}\tau +1$ for each simplex $\tau $ in $K^a.$
\par {\bf 4.2.1. Definitions and Notes.}  A polyhedral spectrum 
$ \{ X_i; f^i_j; {\bf N} \} $ is called standard 
if it satisfies conditions
$S.1, S.3, S.4$ from \cite{koz} and the following condition: 
\par $S.2$ $f^{i+1}_i$ is a simplexial mapping of a complex
$X^a_{i+1}$ onto a normal subdivision (that is a
$p^j$-subdivision with the corresponding $j\in \bf N$)
$X^*_i$ of a complex $X^a_i$ for each $i\in \bf N$ such that 
\par $(i)$ $\lim_{j\to \infty }
\sup _{\tau \subset X^a_j}diam (f^j_i(\tau ))=0$,
where $\tau $ are simplexes in $X^a_j$.
The complex $K^a$ is called uniform if the corresponding to it
polyhedra over $\bf L$ is uniform. 
\par An approximation of all subsets $Y$ of an 
arbitrary ultrauniform space $X$
simultaneously by one spectrum $\{ X_{\alpha }; 
f^{\alpha }_{\beta }; {\sf U} \} $ is called absolute,
that is for each $Y\subset X$ 
there exists ${\sf U}_Y\subset \sf U$ such that
$Y=\lim \{ X_{\alpha }; f^{\alpha }_{\beta }; {\sf U} \} $
and for each $Y'\subset Y$ the following inclusion
${\sf U}_{Y'}\subset {\sf U}_Y$ is satisfied.
This decomposition is called
\par $(a)$ irreducible if for each open $V\subset X$ there exists cofinal
subset ${\sf U}_V\subset \sf U$ such that $ \{ X_{\alpha };
f^{\alpha }_{\beta }; {\sf U}_V \} $ is irreducible.
It is called 
\par $(b)$ uniform if each subspectrum $ \{ X_{\alpha }; 
f^{\alpha }_{\beta }; {\sf U}_V \} $ is uniform.
\par An absolute polyhedral representation is called uniform 
from above, if polyhedra 
$P$ satisfy Condition $2.8.(1.(i))$. 
A polyhedral spectrum is called 
\par $(c)$ non-degenerate if $f^{\alpha }_{\beta }$ are non-degenerate,
that is for each simplex $s$ in $X_{\alpha }$ restrictions
$f^{\alpha }_{\beta }|_s$ are simplexial and 
$dim_{\bf L}f^{\alpha }_{\beta }(s)=dim_{\bf L}s$.
\par Notations and propositions from \S \S II.1-II.10 \cite{koz}
are transferable onto the non-Archimedean case due to
above results and definitions with the substitution of barycentric 
subdivisions onto $p^j$-subdivisions (see for example, 
$\tau ^*_{s+1}$ in \S II.1.11 \cite{koz}).
\par Let $K^c_i$ be defined as in \S II.3.1 \cite{koz}.
For two families $Y^c$ and $Y^c_1$ of simplexes from $K^c_1$
a mapping $h: Y\to Y_1$ is called $c$-linear over 
the field $\bf L$, if for each simplex 
$\tau $ from $Y^c$ (that is from some subcomplex 
$Y^c\cap K^a_{i,i+1}$, where $i\in \bf N$) there exists
$q\in \bf N$ such that $h(\tau )\subset K_{q,q+1}$ and the mapping
$h|_{\tau }$ is simplexial. In view of II.4.7 \cite{koz}
$\{ K_{1,i}; h_{j,i}; {\bf N} \} $ is an inverse spectrum, 
its limit space is denoted $\mbox{ }^{\bf N}K$ and it is called
the cone of the spectrum $\{ K_{1,i}; h_{j,i}; {\bf N} \} $.
Then $\alpha (i)$ is a space of a complete subcomplex
$\alpha ^a(i)\subset X^a_i$, $\alpha ^q(i)=f^{-1}_{i+q,i}(\alpha (i))$,
${\tilde \alpha }^1(i):=X^a_{i+1}\alpha ^1(i)$ is a polyhedral neighbourhood
(clopen star) of a set $\alpha ^1(i)\subset X_{i+1}$ relative
to a complex $X^a_{i+1}$, since the space $X_{i+1}=|X^a_{i+1}|$ 
of a complex $X^a_{i+1}$ is the disjoint union of the corresponding 
simplexes $s$ forming its covering such that spaces $|s|$ 
of these simplexes are clopen in $X_{i+1}$.
Propositions analogous to II.1.9, II.7.18,19, II.9.22 \cite{koz}
follow from Lemma 4.1.9 above.
\par {\bf 4.2.2. Lemma.} {\it $dim_{\bf L}(X_{\alpha }
\cap K_{i-1,i})\le dim_{\bf L}\alpha (i)$.}
\par {\bf Proof.} In view of \S II.7.13 \cite{koz}
we have $X_{\alpha }\cap K_{i-1,i}=X_{\alpha ,i}\cap K_{i-1,i}$.
Let $M^a_{\alpha ,i}$ be a subcomplex of a complex $K^a_{i-1,i}$
such that $M^a_{\alpha ,i}$ is a union of simplexes 
$\sigma \circ f_{i,i-1}(\sigma )$, where $\sigma $ is a simplex
of the complex $X^a_i$ such that $f_{i,i-1}(\sigma )\subset
X_{i-1}\setminus {\tilde \alpha }^1(i-2)$. Then $h_{i,\alpha }:
M^a_{i,\alpha }\to M^a_{i,\alpha }$ is the simplexial mapping
and $h_{i,\alpha }(M_{i,\alpha }\cap \alpha (i))=X_{\alpha ,i}
\cap K_{i-1,i}$ and $h_{i,\alpha }(M_{\alpha ,i}\cap K_{\alpha ,i})
=X_{\alpha ,i}\cap K_{i-1,i}$ (see also \S \S II.9.14-19 \cite{koz}).
For each simplex $s$ we have $dim_{\bf L}s=card (A_s)$, where
$cl( sp_{\bf L}s)=c_0({\bf L},A_s)$, the closure is taken in
the corresponding Banach space $c_0({\bf L},A)$ in which the 
complex is contained, $sp_{\bf L}s$ denotes the $\bf L$-linear 
span of $s$.
\par {\bf 4.2.3. Lemma.} {\it For each $i\in \bf N$ and $\alpha \in \sf U$
there are mappings $h_{i,\alpha }$ satisfying conditions 
\par $(1)$ $h_{i,\alpha }$ is $c$-linear over $\bf L$;
\par $(2)$ $h_{i,\alpha }(x)=h_{i-1,\alpha }(x)$ for each $x\in
K_{1,i-1}$;
\par $(3)$ $h_{i,\alpha }(x)=h_{i-1,\alpha }h_{i,i-1}(x)$ 
for each $x\in h^{-1}_{i,i-1}(\alpha (i-1))$;
\par $(4)$ $h_{i,\alpha }(x)=x$ for each 
$x\in X_i\setminus {\tilde \alpha }^1(i-1)$.}
\par {\bf Proof.}  For $i=1$ the mapping $h_{1,\alpha }$
is identical for each $\alpha \in \sf U$ in accordance with
\S II.9.7 \cite{koz}, hence it satisfies Conditions $(1-4)$.
Let $h_{i,\alpha }$ satisfies Conditions $(1-4)$ for each $i\le n$.
Then $h_{n+1,\alpha }$ is defined by Conditions $(2-4)$.
Moreover, $\alpha (n)$ and
$\tilde \alpha ^1(n)$ may be taken clopen in $X_n$ and $D_1:=h^{-1}_{n+1,n}(
\alpha (n))$ are clopen in $X_{n+1}$, then $D_1$ and $D_2:=K_{1,n}$
and $D_3:=M_{\alpha , n+1}$ are closed in $K_{1, n+1}$. For each simplex
$\tau \subset K_{1, n+1}$ each point $x\in \tau $ has the form
$x=a_1x_1+a_2x_2+a_3x_3$, where $a_j\in B({\bf L},0,1)$, $x_j\in D_j$.
Let $g: c_0\oplus c_0\oplus c_0\to c_0$  be the mapping of evaluation,
that is, $g(x_1,x_2,x_3)=x_1+x_2+x_3$, $h^j: D_j\to K_{1, n+1}$ are 
continuous mappings, so $h_{n+1, \alpha }(x)=(g(a_1h^1(x_1), a_2h^2(x_2),
a_3h^3(x_3))$ is continuous. If to consider uniform polyhedra, then
the simplexial mappings $f: K^a\to N^a$ are 
considered as uniformly continuous.
Then to the cylinder $<N^a,K^a>_f$ it corresponds 
the uniform polyhedron, since
$f(K)=N$. Then $K^a_{i,j}$, $X^a_j$, $X^*_i$, $K^*_{i,j}$ and $K^c_{i,j}$
from \S II.1 in \cite{koz} correspond the uniform polyhedra over $\bf L$,
since $K^*_{i,j}=Sd_{X_j}*K^a_{i,j}$ and $\inf_{(V\cap W=\emptyset , 
V\mbox{ and }W\mbox{ are simplexes from }K^*_{i,j})} dist(V,W)>0$. 
For $T=K_i=\bigcup_{j\in \bf N}K_{j,j+1}$, $K^c_i$, $K^c_{\alpha }$ 
and $X^c_{\alpha }$ the latter condition may be not satisfied, but
$\sup_{(V\mbox{ is simplex from }T)}diam (V)<\infty $. At the same time 
$T$ is not necessarily complete, since it is not necessarily ANRU.
If $f_{j,i}=f^j_i: X_j\to X_i$, $f_j: \mbox{ }^{\bf N}X\to X_i$ are non-stretching
(uniformly continuous), then $h_{j,i}$ from \S II.4 and 
$h_i$ from \S II.5 in \cite{koz} are also 
non-stretching (uniformly continuous
respectively). Then from $h_{1,\alpha }=id$ and $(1-4)$ with
induction by $j$ for each $\alpha $ it follows that $h_{j, \alpha }$
are non-stretching (uniformly continuous correspondingly). 
From this it follows
that the mappings $h_{\alpha }$ from \S II.10 and hence also 
$f_{\alpha , \beta }$ and $f_{\alpha }$ from \S II.11 
and \S II.9.27 \cite{koz} are non-stretching
(uniformly continuous respectively).
\par {\bf 4.2.4. Lemma.} {\it For each open covering $V$ of $W$
which is a clopen subset in a limit 
$\lim_j\{ X_j, f^j_i, {\bf N} \} =:\mbox{ }^{\bf N}X$
of an inverse sequence of polyhedra there are $\alpha \in \sf U$ such that
$W_{\alpha }=W$ and a clopen covering $w_{f_{\alpha }}X^*_{\alpha }
\subset \{ f^{-1}_{\alpha }(X^*_{\alpha }q): q\} $
of $W_{\alpha }$ by $f^{-1}_{\alpha }(X^*_{\alpha }q)$
of the main stars
$X^*_{\alpha }q$ of the complex $X^*_{\alpha }$
is refined into $V$.}
\par {\bf Proof.} Let $\alpha ^a(i)$ be a subcomplex 
of a complex $X^a_i$ such that a simplex $\sigma $
of a complex $X^a_i$ is a simplex of a complex $\alpha ^a(i)$
if and only if $f^{-1}_i(X^a_ie)$ is contained in
the corresponding element $U(e)$ of a covering
$V$ for each its vertex $e$. Then $\alpha ^a(i)$ is a 
complete subcomplex of a complex $X^a_i$. That for a sequence
$\alpha (i)$ to verify condition 
\par $(A)$ $\mbox{ }{\tilde \alpha }^1(i)\subset \alpha (i+1)$
for each $i\in \bf N$ as 
in \S II.11.20 \cite{koz}
it is sufficient to mention that $X^*_i$ is the $p^j$-subdivision 
of the complex
$X^a_i$ with $j\ge 1$, so $X^a_i\tau \subset X^a_iq$ for some vertex
$q\in X^a_i$. In this case the main stars are clopen, hence 
$w_{f_{\alpha }}X^*_{\alpha }$ is the clopen covering.
If $ \{ q_j:$ $j=0,1,2,... \} $ are vertexes 
from $X^a_i$ stars of which contain
$f_i(x)$, then there is a clopen simplex $\tau \subset X^a_i$ with such
vertexes and $f_i(x)\in \tau $. At the same time the star $St(q_j,
w_{f_i}X^a_i)$ with each $q_j$ is contained in $V$, consequently,
$\tau \subset \alpha ^a(i)$ and $f_i^{-1}(\alpha (i))$ is clopen
in $W_{\alpha }$. 
Each $p^j$-subdivision $X^*_{i-1}$ of a complex $X^a_{i-1}$
is proper for each $j\in \bf N$, consequently,
$X^*_{i-1}\sigma \subset X^a_{i-1}q$ for some suitable vertex $q$ in 
the complex $X^a_{i-1}$, where $\sigma $ is the simplex in
$X^*_{i-1}$. Therefore, $h_{i-1}(X^*_{\alpha }e)\subset X^*_{i-1}q$,
where $q$ is a vertex of the complex $\alpha ^a(i-1)$.
Since $f_{\alpha }^{-1}(X^*_{\alpha }q)$ contains
the clopen subcovering, then $W_{\alpha }$ is clopen and
$W_{\alpha }=W$. In the case of separable $f_{\alpha }^{-1}(X^*_{\alpha }q)$
this subcovering can be chosen disjoint due to ultrametrizability
of $f_{\alpha }^{-1}(X^*_{\alpha }q)$ (see Theorem 7.3.3 in \cite{eng}).
\par Vice versa, let $x\in W$ and $x\in U\in V$. Coverings
$\omega _{f_i}X^*_{\alpha }$ satisfy Condition 4.2.1.(i) such that 
$St(x,\omega _{f_i}X^a_i)\subset U$. Let $\{ e_j:$ $j\in \Lambda _{i,x} \} $
be vertexes of the complex $X^a_i$, stars of which contain 
the point $f_i(x)$. Then $f_i(x)\in \tau $, where $\tau $
is a clopen simplex of a complex $X^a_i$ with the described above vertexes.
Since $f_i^{-1}(X^a_ie_j)\subset U$ for each $j\in \Lambda _{i,x}$,
then $\tau $ is the simplex in the complex $\alpha ^a(i)$.
Therefore, $x\in f_i^{-1}(\tau )\subset f_i^{-1}(\alpha (i))
\subset V_{\alpha },$ hence $W=V_{\alpha }$.
\par {\bf 4.2.5. Lemma.} {\it If $Y\subset \mbox{ }^{\bf N}X$, 
$f_i(Y)=X$, $f^j_i:
X_j\to X_i$ are irreducible mappings for each $j\ge i\in \bf N$, then 
$f^{\beta }_{\alpha }: X^*_{\beta }\to X^*_{\alpha }$ are irreducible
and $f_{\alpha }(Y\cap V_{\alpha })=X^*_{\alpha }$ for each $\alpha \in
\sf U$ and $Y$ is dense in $\mbox{ }^{\bf N}X$. If $f_i$ are irreducible
for each $i$, then $f_{\alpha }$ are irreducible for each 
$\alpha \in \sf U$.}
\par {\bf Proof.} From irredicibility of $f^j_i$ and 
Lemma 3.5 it follows that
for each (clopen) simplex $\tau \subset X^a_i$ the mapping $f^j_i|_{
(f^j_i)^{-1}(\tau )\subset X_j} $ is essential. Let $\beta \in \sf U$
and $t$ be a simplex in the complex $X^*_{\beta }$, then there exists
some $p^j$-subdivision $\sigma ^*$ of the simplex $\sigma $ from 
$K^a_{i-1,i}\cap X^c_{\beta , i}$ and there exists a simplex $\tau $ from
$M_{\alpha ,i}\cap \alpha (i)$ such that $\sigma =h_{i, \alpha }(\tau )$
and $t$ is a simplex from $\sigma ^*$. From \S II.3.2,10.2 ($c$-linearity
of $h_{\alpha }: K_1\to K_1$), \S II.7.1.4,7.15 and \S II.4 ($f_{\beta ,
\alpha }=h_{\alpha }|_{X_{\beta }}$, $X_{\alpha }=\bigcup_{i\in \bf N}
X_{\alpha ,i}$, $X^c_{\alpha ,i}=X^c_{\alpha }\cap K^c_{1,i}$) it follows
that $f^{\beta }_{\alpha }=f_{\beta , \alpha }$ is irreducible due to Lemmas
3.4 and 3.5, since $h_{\alpha }$ is essential on $t\subset X_{\beta }$.
Indeed, due to II.4.9 we have $h_{j,i}=f_{j,i}$ for $x\in X_j$ and $j\ge i$.
In view of II.10.3 we have $h_{\alpha }(x)=h_{i,\alpha }(x)$ for 
$x\in K_{1, \alpha }$. In view of II.9.1 we have $h_{\alpha }(x): K_{1,i}
\to K_{1,i}$ is $c$-linear and $h_{1,\alpha }=id$, in addition, by
II.9.3 $h_{i,\alpha }(x)=h_{i-1, \alpha }\circ h_{i,i-1}(x)$ for $x\in
h^{-1}_{i,i-1}(\alpha (i-1))$, by II.9.4 $h_{i, \alpha }(x)=x$ for 
$x\in X_i\setminus \tilde \alpha ^1(i-1)$, also $\tilde \alpha ^1
(i)=X^a_{i+1}\alpha ^1(i)$ is the clopen star relative to the complex
$X^a_{i+1}$, where $\alpha ^1(i)$ is the clopen subcomplex 
(and the subpolyhedron as the subspace) in $X^a_{i+1}$, consequently,
we can suppose that $\tilde \alpha ^1(i)=\alpha ^1(i)$ 
in this non-Archimedean
case. The rest of the statements of this Lemma follows from 
\S II.11.22 in \cite{koz}.
\par {\bf 4.2.6. Definition.} A mapping $f: X\to K^a$ into a 
simplexial complex $K^a$ is called locally finite, if 
each $x\in X$ has a neighbourhood $U$ such that
$f(U)$ is contained in a finite subcomplex 
(that is the subcomplex consisting of finite number of
simplexes and their verges) of a complex $K^a$.
\par {\bf 4.2.7. Lemma.} {\it If all mappings $f_i: 
\mbox{ }^{\bf N}X\to X_i$ are locally finite and $X_i$
are locally finite-dimensional over $\bf L$, then
all mappings $f_{\alpha }: V_{\alpha }\to X^*_{\alpha }$
are locally finite.}
\par {\bf Proof.} Since $X_i$ are locally finite-dimensional
over $\bf L$, then each its $p^j$-subdivision
is also locally finite-dimensional over $\bf L$, hence 
$X^*_i$ and $X^c_{\alpha }$ are locally finite-dimensional
over $\bf L$, since $\bf L$ is locally compact. 
For $x\in V_{\alpha }$ and $\alpha \in \sf U$
we have $x\in f^{-1}({\tilde \alpha }^1(i))$, since
$V_{\alpha }=\bigcup_{i=1}^{\infty }f_{i+1}^{-1}({\tilde \alpha }^1(i))$.
If $W$ is a neighbourhood of $x$ such that 
$f_{i+1}(W)$ is a finite subcomplex $P^a$ of a complex
$X^a_{i+1}$. Taking $W\subset f_{i+1}^{-1}(\alpha (i+1))$, 
then $P\subset \alpha (i+1)$. Since $h_{\alpha }$ is
$c$-linear over $\bf L$, then $h_{\alpha }(P)$ 
is contained in a finite number of simplexes from
$X^c_{\alpha }$, hence in a finite number of simplexes of its
subdivision $X^*_{\alpha }$, hence $f_{\alpha }$ is locally finite, since
$f_{\alpha }|_{f_{i+1}^{-1}(\alpha (i+1))}=h_{\alpha }f_{i+1}|_{
f_{i+1}^{-1}(\alpha (i+1))}$ due to Propositions 
II.10.4 and II.11.3 \cite{koz}.
\par {\bf 4.2.8. Lemma.} {\it If all $X_i$ are finite-dimensional
$dim(X_i)<\infty $, then 
the complexes $X_{\alpha }$ can be chosen locally finite and 
such that $dim_{\bf L}
(X_{\alpha })\le \sup_{i\in \bf N}dim_{\bf L}X_i$.}
\par {\bf Proof.} This follows from Theorems 3.1, 3.18,
since complexes $K_1$ for each $Y_j$ can be chosen locally 
finite and $X_{\alpha }\subset K_1$. The second statement follows from
Lemma 4.2.2.
\par {\bf 4.2.9. Theorem.} {\it Each ultrauniform space $(X, {\bf P})$ 
has an irreducible absolute normal and uniform from above polyhedral 
representation 
$$T:= \{ X_{\beta }, f^{\alpha }_{\beta }, {\sf U} \} $$
over $\bf L$, that is, for each $Y\subset (\tilde X, {\bf P})$
there exists an ordered subset ${\sf U}_Y\subset \sf U$ such that
there exists an uniform isomorphism 
$$\mbox{ }^Yf: Y\to
\lim \{ X_{\beta }, f^{\beta }_{\alpha }, {\sf U}_Y \} $$ 
(that is $\mbox{ }^Yf$ is a bijective surjective mapping, $\mbox{ }^Yf$
and $\mbox{ }^Yf^{-1}$ are uniformly continuous).
Moreover,
as an absolute polyhedral expansion $T$ (not necessarily uniform
such that $X$ is homeomorphic with $\lim T$) may be 
chosen locally finite-dimensional $X_{\beta }$ over the field $\bf L$.
If $dim(X)=k$, then there exists $k$-dimensional
$T$ over $\bf L$ and 
$$ \{ X_{\beta }, f^{\alpha }_{\beta }, {\sf U}_X \} $$
is non-degenerate (and irreducible for complete $X$) such that
there exists a homeomorphism 
$$g: X\to \lim 
\{ X_{\beta }, f^{\alpha }_{\beta }, {\sf U}_X \} $$
with irreducible surjective mappings $g_{\alpha }:
=f_{\alpha }\circ g$.}
\par {\bf Proof.} Let $S=\{ X_j, f^j_i \}$ be an inverse polyhedral sequence
over $\bf L$. Then sections $X_{\alpha }$ and bonding mappings 
$f^{\beta }_{\alpha }$
for $\beta \ge \alpha \in \sf U$ form the inverse mapping system 
$$R=\{ X_{\beta }, f^{\beta }_{\alpha }, {\sf U} \} $$ 
that is an expansion of
$S$ and the absolute polyhedral expansion for $X=\lim S$.
If $X_i$ are uniform polyhedra and $f^j_i$ are non-stretching 
(or uniformly continuous) irreducible simplexial mappings for each $j\ge i
\in \bf N$, then the same is true for $f^{\beta }_{\alpha }$ for 
each $\beta \ge \alpha \in \sf U$. From Theorem 3.18, Lemmas and Notes
4.1.1-4.2.8 
above and \S II.12 in \cite{koz} it follows the existence of the 
homeomorphism 
$$\mbox{ }^Yf: \lim \{ X_{\beta }, f^{\beta }_{\alpha } ,{\sf U}_Y \} \to Y$$ 
such that $\mbox{ }^Yf$ 
is uniformly continuous, since $f^{\beta }_{\alpha }$ 
and $f_{\alpha }$ are non-stretching. From Lemma 3.17 it follows
that $\mbox{ }^Yf$ is the uniform embedding, since 
$$ \lim_{i\to \infty }
\sup_{(s\mbox{ are simplexes from }X^a_i)}diam(s)=0$$ 
and due to the consideration
of the sequence $\{ \alpha (j): j\ge i \} =\alpha $. 
Indeed, the corresponding
$X_{\alpha }$ has the cofinal subset $E\subset {\sf U}_Y$ such that for each
$b>0$ the set $\{ \alpha : \sup_{(s\mbox{ are simplexes from }X_{\alpha })}
diam(s)>b$ is finite and $V_{\alpha }=h^{-1}_{\alpha }(X_{\alpha })\cap
\mbox{ }^{\bf N}X$ and $Y=\bigcap_{\alpha \in {\sf U}_Y}V_{\alpha }$,
so that $\mbox{ }^Yf^{-1}$ is uniformly continuous.
\par Each ultrametric space $(X, \rho )$ has the topological embedding
into $B({\bf L},0,1)^{w(X)}$ and into $D^{w(X)}$ which are compact (see
\S 6.2.16 in \cite{eng}). In particular, the separable space $X$ has
the embedding into $B({\bf L},0,1)^q$ for each $q\in \bf N$ and 
for $q=\aleph _0$, but not uniformly in general, for example, 
$c_0({\bf L},{\bf N})$ is complete, but it is not compact. From this
the last statement of the Theorem follows. That is, for 
the polyhedral expansions reproducing the uniformity up to 
an isomorphism
with the initial ultrauniform space in general infinite-dimensioanl over
$\bf L$ polyhedra may be necessary and $\sup_{\beta }dim_{\bf L}X_{\beta }
\le w(X)$. For the polyhedral expansion reproducing $X$ up to 
a (topological) homeomorphism
locally finite-dimensional $X_{\beta }$ are sufficient, 
that is for each $x\in X_{\beta }$ there exists a clopen simplex $s$ in 
$X_{\beta }$ such that $x\in s$ and $dim_{\bf L}s<\aleph _0$,
since ${\bf L}^{w(X)}$
is a limit of an inverse system of finite-dimensional polyhedra
over $\bf L$.
In view of Lemma 2.5 
it is sufficient to take $\sf U$ with $card ({\sf U})\le 
\aleph _0card ({\bf P'})2^{w(X)}$, where $\bf P'$ is a subfamily of 
$\bf P$ of minimal cardinality, which generate the same uniformity
in $X$, since sections $X_{\alpha }$ of the spectrum for each 
ultrametric space $Y_j$ are spaces of $X^c_{\alpha }$
for $Y_j$, where $X^c_{\alpha }$ is a collection
of simplexes $s$ from $K^c_1=\bigcup_{q=1}^{\infty } K^a_{q,q+1}$
such that all vertexes of $s$ are in $\bigcup_{q=1}^{\infty }(\alpha (q)
\setminus {\tilde \alpha }^1(q-1))$.
If $card (X)\ge \aleph _0$, then $\aleph _0card ({\bf P'})
2^{w(X)}=2^{w(X)}$, since $2^{w(X)}>\aleph _0 card ({\bf P'})$.
\par If $dim(X)=k$, then by Theorem 3.18 there exists $k$-dimensional 
over $\bf L$ inverse mapping system $S$ associated with coverings
of order $k$. 
For $dim(X)=\infty $ we can use orders of coverings
as cardinals equal to $k\ge \aleph _0$ and as follows from
results above bonding mappings $f^n_m$ and polyhedra $P_m$
can be chosen such that there exists $m\in {\sf U}_X$ with
$dim_{\bf L}P_m=k$. The final part of the proof is analogous to that of
Theorem 12.21 \cite{koz}.
\par {\bf 4.2.10. Corollary.} {\it If in a complete ultrauniform space 
$(X,{\bf P})$ a subspace $R$ is compact and each closed subset $G$, $G\subset
X\setminus R$ is such that $G$ may be uniformly embedded into
$B({\bf L^n},0,1)$ with $n=n(G)\in \bf N$, then $X$ is uniformly isomorphic 
with $\lim S$, where $S=\{ P_m, f^m_n, E \}$ is an irreducible normal 
inverse mapping system and $P_m$ are uniform finite-dimensional over $\bf L$
polyhedra.}
\par {\bf Proof.} For $X_{\rho }\setminus R_{\rho }$ there is a sequence
of closed subsets $G_{\rho }(k)$ with $\bigcup_{k\in \bf N}G_{\rho }(k)=
X_{\rho }\setminus R_{\rho }$ and uniform clopen neighbourhoods
$G_{\rho }(k)\subset U_{\rho }(k)\subset B({\bf L^n},0,r_n)$, $n=n(G_{\rho }
(k))$ such that there exists a uniformly continuous retraction $r_{k, \rho }:
U_{\rho }(k)\to G_{\rho }(k)$ and a family 
$$b_{\rho ,k}:=
\sup_{y\in U_{\rho }(k)}\inf_{x\in G_{\rho }(k)}\rho (x,y)$$
satisfying conditions of Lemma 3.17. From $R_{\rho }\subset 
B({\bf L^n},0,1)$ and the existence of the family $\{ U_{\rho }(k) \}$
such that $\lim \{ U_{\rho }(k); f^{\rho ,k}_{\rho ',k'}, F \}$ 
$=X\setminus Int(R)$ for some directed $F\subset {\bf P}\times {\bf N}$
it follows the existence of the inverse spectra $S$. Indeed, $f_{\rho ,k}$
may be extended on $R\setminus Int(R)$ and compacts 
$f_{\rho ,k}(R\setminus Int(R))$ may be uniformly 
embedded into $U_{\rho }(k)$,
and spectra for $R$ and $X\setminus Int(R)$ may be done consistent on
$R\setminus Int(R)$ that to get $S$ for $X$
(see also Lemmas 3.2, 3.3.1).
\par {\bf 4.2.11. Note.} In view of Theorem 1.7 in \cite{isb2} an incomplete 
uniform space $(X,{\bf P})$ is not ANRU, hence in $S$ in general polyhedra
$X_{\alpha }$ are incomplete, for example, $X_{\alpha }={\bf Q}\times
B({\bf L},0,1)$ in $\bf L^2$.
\par {\bf 4.3. Theorem.} {\it Let $G$ be a locally compact group with a 
left ultrauniformity $\bf P$ and $X$ is closed in $G$. Then $X$ is
uniformly isomorphic with a limit of an inverse mapping system
$\lim S$ and $P_m$ are
finite-dimensional uniform polyhedra over $\bf L$ of an irreducible normal 
inverse mapping system $S=\{ P_m, f^m_n, E \} $.}
\par {\bf Proof.} Let $V$ be a compact clopen neighbourhood of the unit 
element $e\in G$, $\bf P$ be the family of the left -invariant 
ultrapseudometrics in $G$ ($\rho (yx,yz)=\rho (x,z)$ for each
$x,$ $y,$ $z\in G$ and $\rho \in \bf P$,
see \cite{eng,itrsw}). That is $V=\{ x\in G:
\rho (x,e)\le b \} $ for some $\rho \in \bf P$ and $b\in \Gamma _{\bf L}$.
The closed subgroup $H_{\rho }:=\{ x\in G: \rho (x,e)=0 \}$ defines the 
equivalence relation $T_{\rho }$ in $G$: $gT_{\rho }h$, if $h^{-1}g
\in H_{\rho }$, since the left classes $hH_{\rho }$ are either disjoint
or coinside. Therefore, there exists the uniformly continuous quotient 
mapping $\pi _{\rho }: G\to G_{\rho }:=G/T_{\rho }$, in addition
$V\supset H_{\rho }$ and $VH_{\rho }=V$. From the compactness of 
$V/T_{\rho }=:V_{\rho }$ it follows that $G_{\rho }$ is a locally
compact ultrametric space. From the completeness of $X$ it follows
that $X_{\rho }$ is closed in $G_{\rho }\hookrightarrow 
c_0({\bf L},A_{\rho })$ with $card(A_{\rho })=w(G_{\rho })$. Therefore,
$X\cap V$ has due to Theorem 3.18 the finite-dimensional over $\bf L$
polyhedral expansion. In $G$ there is a family $M$ of elements
such that $gV\cap hV=\emptyset $ for each $g\ne h\in M$ and 
$G=\bigcup_{g\in M}gV$. If $g\ne h\in M$, then $\rho (g,h)>b$, since
from $\rho (g,h)\le b$ it would follow that $\rho (gV,hV)\le b$ and 
$gV\cap hV\ne \emptyset $, that contradicts the definition of $M$.
Then $Y(g)\cap Y(h)=\emptyset $ for $g\ne h\in M$ with $Y(g):=
\pi _{\rho }(gV)$ and $G_{\rho }=\bigoplus_{g\in M}Y(g)$ and $X_{\rho }
=\bigoplus_{g\in M} X_{\rho }\cap Y(g)$. From this the statement of 
the Theorem follows.
\par {\bf 4.4. Theorem.} {\it An ultrauniform space $(X,{\bf P})$ 
is homeomorphic with $\lim S$ with \v Cech-complete finite-dimensional
over $\bf L$ polyhedra $P_m$ and an irreducible inverse mapping system
$S=\{ P_m, f^m_n, F \}$ if and only if $X$ is \v Cech-complete.}
\par {\bf Proof.} At first we verify that for $(X,{\bf P})$ the
condition of \v Cech-completeness is equivalent with
the $\bf L$-completeness.
The latter means that $X$ has the embedding
as the closed subspace into $\bf L^{\tau }$,
where $card(\tau )\le w(X)$. If $X$ is $\bf L$-complete, then $X$
has the embedding as the closed subspace into $\bf R^{\tau }$, 
since $\bf L$ is
Lindel\"of and $X$ is \v Cech-complete due to \S 3.11.3, 3.11.6 and 3.11.12
in \cite{eng}.
\par To prove the reverse statement we verify that the $\bf L$-completeness
is equivalent to the satisfaction of the following conditions:
there is not any Tychonoff space $Y$ for which 
\par $(LC1)$ there is the 
homeomorphic embedding $r: X\to Y$ with $r(X)\ne cl(r(X))=Y$ and 
\par $(LC2)$
for each continuous function $f: X\to \bf L$ there is a continuous
function $g: Y\to \bf L$ such that $g\circ r=f$. 
Evidently, that it is sufficient
to consider the case of ultrametrizable $X$. For such $X$ the family of
continuous $f: X\to \bf L$ separates points $x\in X$ and closed subsets
$C$ for which $x\notin C$, for example, $f(x)=\rho (x,C)=\inf_{y\in C}
\rho (x,y)$. Then substituting $\bf R$ on $\bf L$ in the proof of Theorem
3.11.3 \cite{eng} we get the desired statement.
\par There is a continuous function $h: {\bf L}\to \bf R$ with
$h({\bf L})=\bf R$ and for each continuous $q: X\to \bf R$ there is a 
continuous function $t: X\to \bf L$ such that $h\circ t=q$. The latter
is sufficient to verify for $(X,\rho )$, since for $(X, {\bf P})$
this statement may be obtained with the limits of the inverse system.
For $(X, \rho )$ and each $b_n>0$ there exists a locally constant
$t_n: X\to \bf L$ such that $h\circ t_n$ and $q$ are $b_n$-close, since
$(X, \rho )$ is the disjoint union of clopen subsets $Y(j)$ with
$\sup_j diam (f(Y(j))\le b$, $h^{-1}(Y(j))\subset \bf L$, $h^{-1}(Y(j))
\cap h^{-1}(Y(i))=\emptyset $ for $i\ne j$. Therefore, there exists
a sequence $\{ t_n: n\in {\bf N} \} $ giving the desired 
$t=\lim_{n\to \infty }t_n$, where $\lim_{n\to \infty }b_n=0$.
Let $(X, {\bf P})$ is \v Chech-complete and not $\bf L$-complete, 
then there would be a homeomorphic embedding $r: X\to Y$ with
$r(X)\ne cl(r(X))=Y$ and for each continuous function $f: X\to \bf L$
there is a continuous function $g: Y\to \bf L$ with $g\circ r=f$. 
Consequently, each continuous function $q: X\to \bf R$ has continuous
functions $t: X\to \bf L$ and $g: Y\to \bf L$ with $g\circ r=t$ 
and $h\circ g\circ r=h\circ t=q$, that contradicts \v Chech-completeness
of $X$.
\par If $X$ is homeomorphic with $\lim S$, then $X$ is \v Chech-complete,
since each $P_m$ is Lindel\"of and \v Chech-complete. 
\par If $X$ is \v Chech-complete, then $X$ is homeomorphic with
$\lim \{ X_{\alpha }, g^{\alpha }_{\beta }, E\}$, where each 
$X_{\alpha }$ is locally compact and has the embedding into 
$\bf L^{n(\alpha )}$ as the clopen subspace, where $n(\alpha )\in \bf N$
(see also Corollary 5 from Theorem 9 in \cite{pas} and Proposition 
3.2.17 in \cite{fed}). From Theorems 3.18 and 4.9 it follows that
$X$ is homeomorphic with $\lim S$, where $P_m$ are finite-dimensional
over $\bf L$ polyhedra.
\section{Polyhedral expansions and 
relations between ultrauniform and uniform spaces.}
\par {\bf 5.1. Theorem.} {\it For each ultrauniform space $(X,{\sf U})$
there exists a uniform locally connected 
space $(Y,{\sf V})$ and a continuous
quotient mapping $\theta : X\to Y$ such that 
\par $(i)$ $\theta (X)=Y$, $d(X)=d(Y)$, $w(X)=w(Y)$; 
\par $(ii)$ if $X$ is dense
in itself, then $Y$ is dense in itself;
\par $(iii)$  if $X=\bigoplus_jX_j$,
then $Y=\bigoplus_jY_j$;
\par $(iv)$ there exist natural compactifications
$cX$ and $c_{\theta }Y$ such that $\theta $ has a continuous extension
on $cX$ and $\theta (cX)=c_{\theta }Y$;
\par $(v)$ if $Z$ is a subspace of $X$, then there exists a subspace
$W$ of $Y$ such that Condition $(i)$ is satisfied for $Z$ and $W$
instead of $X$ and $Y$;
\par $(vi)$ if $X$ has a $n$-dimensional polyhedral expansion
over $\bf L$, then $Y$ has a $n$-dimensional 
polyhedral expansion over $\bf R$.}
\par {\bf Proof.} In view of Theorem 2.3.24 \cite{eng} the Cantor cube
$D^{\sf m}$ has the weight $\sf m$, where $D$ is the two-point 
discrete set, $D^{\sf m}=\prod_{s\in S}D_s$, $card (S)=\sf m$, 
$D_s=D$ for each $s\in S$. In view of Corollary 6.2 \cite{eng} for each 
non-Archimedean infinite locally compact field $\bf L$ 
with a non-trivial valuation a ball $B({\bf L},x,r):=\{ y\in {\bf L}:
|x-y|\le r \} $ is homeomorphic to $D^{\aleph _0}$. On the other hand, there
exists a well-known quotient mapping of $D^{\aleph _0}$ onto $[0,1]$. It
can be also realized for $B({\bf L},0,1)$
and in particular for ${\bf Z_p}:=B({\bf Q_p},0,1)$, where $p$
is a prime. Let $s$ be a unit ball in $c_0({\bf L},\alpha )$ and
$\theta : B({\bf L},0,1)\to [0,1]$ be a quotient mapping. 
Then this mapping generates
a mapping $\theta : s\to \theta (s)\subset c_0({\bf R},\alpha )$, where
$\theta (s)=\{ y\in \theta (x):$ $x=(x^j: j\in\alpha ), |x^j|\le 1,
x\in c_0({\bf L},\alpha ) \} $, hence $(y^j: j\in \alpha )=y$,
$y^j=\theta (x^j)$ and for each $\epsilon >0$ we have $card \{ j:
|x^j|>\epsilon \} <\aleph _0$, hence $card \{ j: |y^j|>\epsilon \} <\aleph _0$
such that $y^j\in [0,1]$. If $\{ e_j: j\in \alpha \} $ is the standard 
orthonormal basis in $c_0({\bf R},\alpha )$, then under the embedding 
$c_0({\bf R},\alpha )\hookrightarrow c_0({\bf R},\alpha )^*$ associated 
with this basis to each $e_j$ there corresponds a continuous linear functional
$e^j\in c_0({\bf R},\alpha )^*$, where $c_0({\bf R},\alpha )^*$
is a topological dual space, such that $e^j(e_i)=\delta ^j_i$.
Therefore, $\theta (s)=\{ y\in c_0({\bf R},\alpha ):$ $0\le e^j(y)\le 1 \} $, 
which is by definition a simplex in $c_0({\bf R},\alpha )$ 
such that if $y\in \theta (s)$ and $0< \| y \| <1$, then
$B(c_0({\bf R},\alpha ), y, \min (\| y\| ,1- \| y \| ))\subset
\theta (s) .$
\par On the other hand, if $b$ is a simplex in $c_0({\bf R},\alpha )$
and $b\subset B(c_0({\bf R},\alpha ),0,1)$, then $\theta ^{-1}(Int(b))$
is open in $c_0({\bf R},\alpha )$, $\theta ^{-1}(Int(b))\ne \emptyset $
if $Int(b)\ne \emptyset $, since $\theta $ is continuous.
Choose $\alpha $ such that $b$ is a canonical closed subset in
$c_0({\bf R},\alpha )$, then $\theta ^{-1}(b)$ is closed in
$c_0({\bf L},\alpha )$. Therefore, $cl(\theta ^{-1}(Int(b)))=\theta ^{-1}(b)$,
since $\theta ^{-1}(Int(b))\subset \theta ^{-1}(b)$ and if $c$ is a canonical
closed subset in $[0,1]$, then $\theta ^{-1}(c)$ is a canonical closed subset
in $D^{\aleph _0}$.
\par If $X=\lim \{ X_{\alpha } ,f^{\alpha }_{\beta },E \} $ is a spectral 
decomposition of an ultrauniform space and $X_{\alpha }$ are polyhedra in
$c_0({\bf L},\gamma _{\alpha })$, 
then $\theta $ generates $\tilde \theta _{\alpha }
: X_{\alpha }\to \tilde \theta _{\alpha }(X_{\alpha }),$ where
$\tilde \theta _{\alpha }(X_{\alpha })$ are polyhedra in
$c_0({\bf R},\gamma _{\alpha })$, since each $X_{\alpha }$
can be presented as a union of simplexes of diameter not greater, than $1$.
Then ${\tilde f}^{\alpha }_{\beta }:={\tilde 
\theta }_{\beta }\circ f^{\alpha }_{\beta }\circ 
{\tilde \theta }_{\alpha }^{-1}:$ ${\tilde \theta }_{\alpha }(X_{\alpha })
\to {\tilde \theta }_{\beta }(X_{\beta })=:{\tilde X}_{\beta }$
are continuous bonding mappings. In view of Lemma 2.5.9 \cite{eng}
the topological space
${\tilde \theta }(X)=\lim \{ {\tilde X}_{\alpha },
{\tilde f}^{\alpha }_{\beta }, E \} $
is locally connected and ${\tilde \theta }(X)$ is the quotient space
of $X$. 
\par Now we proof that if $X$ is given, then 
$Y$ satisfies Conditions $(i-vi)$.
In view of Proposition 2.5.6 \cite{eng}
are satisfied $(iii,v)$. On the other hand,
$dim_{\bf L}X_{\alpha }=dim_{\bf R}{\tilde X}_{\alpha }$
for each $\alpha \in E$, consequently, $(vi)$ is fulfilled.
Since $d(X_{\alpha })=d({\tilde X}_{\alpha })$
and $w(X_{\alpha })=w({\tilde X}_{\alpha })$
for each $\alpha $, then $(i)$ is fulfilled.
If $x$ is an isolated point in $X$, then
there exists $\alpha \in E$ such that $x_{\alpha }$ is an isolated
point of $X_{\alpha }$, consequently, 
${\tilde x}_{\alpha }$ is an isolated point of 
${\tilde X}_{\alpha }$.
\par For prooving $(iv)$ we consider 
$cX$ equal to some specific compactification
${\tilde c}{\tilde X}$ of $\tilde X$, where $\tilde X$ 
is a completion of the uniform space $X$.
Therefore, the consideration reduces to the case of a complete
$X$. In such case $X$ has the uniform polyhedral expansion
such that conditions $2.8(i,ii)$ are satisfied.
Then we take $cX$ 
as $\lim \{ c_{\alpha }X_{\alpha },{\bar f}^{\alpha }_{\beta },E \} $,
where $c_{\alpha }X_{\alpha }$ are compact spaces.
Since each $X_{\alpha }$ is a uniform polyhedron
we take $c_{\alpha }X_{\alpha }$ as a topological sum
of compactifications of disjoint polyhedra $s$ in $X_{\alpha }$.
Practically we take the Banaschewski compactification
of each $s$ \cite{roo}. If $s=B(c_0({\bf L},\gamma _{\alpha }),x,r)$,
then its Banaschewski compactification $s^{\zeta }$ is isomorphic with
$B({\bf L},0,1)^a$, where $a=2^{\gamma _{\alpha }\aleph _0}$.
Therefore, ${\tilde \theta }_{\alpha }$
has a continuous extension onto $c_{\alpha }X_{\alpha }
=\bigoplus_{s\in X_{\alpha }}c_{\alpha }s$.
Evidently, ${\tilde \theta }_{\alpha }(s)=[0,1]^a$
is the Stone-\v Cech compactification of ${\tilde \theta }(s)$.
\par It remains to prove statements about compactifications of simplexes.
Let $Bco(s)$ denotes an algebra of clopen subsets of
$s$. Instead of a ring homomorphism $\eta : Bco(s)\to \bf F_2$
in the Banaschewski construcion let us consider a ring homomorphism
$\eta : Bco(s)\to \bf Z_p$, where $\bf F_2$ is a finite field consisting
of two elements. Topologically $\bf Z_p$ is homeomorphic to
${\bf F_2}^{\aleph _0}$, where $a={\sf c}^{w(s)}=
2^{\gamma _{\alpha }\aleph _0}$. 
The family $\sf \Omega $ of all such distinct
$\eta $, which has $\sigma $-additive extensions onto 
$Bf(s)$ has the cardinality $|{ \sf \Omega  } |={\sf c}^{w(s)}$,
where $w(s)=\gamma _{\alpha }\aleph _0$
and $Bf(s)$ denotes the Borel $\sigma $-field of $s$, 
since $|{\bf Z_p}|=\sf c$ and each such $\eta $ is characterised
by a countable family of atoms $a_j$ such that $\lim_{j\to \infty }
\eta (a_j)=0$ (see Ch. 7 \cite{roo}).
Therefore, $s$ has an embedding into ${\bf Z_p}^a={\bf F_2}^a$.
The family ${\sf \Omega }$ is a subfamily of 
a family $\sf \Upsilon $ of all ring homomorphisms $\eta $.
Another way to realize a totally disconnected compactification
$\gamma s$ of $s$ is in the consideration
of the space $H:=C^0_b(s,{\bf Q_p})$ of bounded continuous functions
$f: s\to \bf Q_p$. To each $x\in s$ there corresponds 
a continuous linear functional $x'\in H^*$ such that
$x'(f)=f(x)$ and $\| x' \|=1$ due to Hahn-Banach Theorem,
since $\bf Q_p$ is spherically complete, where $H^*$ denotes 
the Banach space of continuous linear functionals on $H$. In view of 
Alaoglu-Bourbaki Theorem (see Exer. 9.202 \cite{nari})
the unit ball $U:=B(H^*,0,1)$ is $\sigma (H^*,H)$-compact.
Therefore, the closure $\gamma s$ of $\gamma (s)$ in $U$ is compact,
where $\gamma $ is an embedding of $s$ into $U$. From this construction
we have $\gamma s={\bf Z_p}^a$.
From Example 9.3.3 \cite{nari} it follows that the 
Stone-\v Chech compactification of $\theta (s)$ is
equal to $[0,1]^a$. Then $s^{\zeta }=\gamma s$, since
$d(s^{\zeta })=d(s)$ and due to Hewitt-Marczewski-Pondiczeri 
Theorem \cite{eng} ${\bf Z_p^a}\subset s^{\zeta }\subset {\bf Z_p^a}$.
\par {\bf 5.2. Theorem.} {\it For each locally connected
uniform space $(Y,{\sf V})$ there exists an ultrauniform
space $(X,{\sf U})$ and a continuous
quotient mapping $\theta : X\to Y$ such that 
\par $(i)$ $\theta (X)=Y$, $d(X)=d(Y)$, $w(X)=w(Y)$; 
\par $(ii)$ if $Y$ is dense
in itself, then $X$ is dense in itself;
\par $(iii)$  if $Y=\bigoplus_jY_j$,
then $X=\bigoplus_jX_j$;
\par $(iv)$ there exist natural compactifications
$cX$ and $c_{\theta }Y$ such that $\theta $ has a continuous extension
on $cX$ and $\theta (cX)=c_{\theta }Y$;
\par $(v)$ if $Y$ has a $n$-dimensional polyhedral expansion
over $\bf R$, then $X$ has a $n$-dimensional 
polyhedral expansion over $\bf L$.}
\par {\bf Proof.} From the proof of Theorems 2.4 and 5.1
it follows that there exists a topological 
embedding of $Y$ into a subspace $Z$ of $H^*$, where
$H:={\tilde C}^0_b(Y,{\bf R})$, such that $Z$ is isomorphic with
$c_0({\bf R},ord(Y))$, where $ord(Y)$ is an ordinal
associated with $Y$, $Z=cl(sp_{\bf R}Y)$, the closure of the 
linear span of $Y$ over $\bf R$ is taken in $H^*$
(see also \cite{eng,isb3,lu3}).
Therefore, polyhedra of a polyhedral expansion of $Y$
can be realized in $Z$. To each simplex
$q$ in $c_0({\bf R},ord(Y))$ there corresponds a simplex $s$
in $c_0({\bf Q_p},ord(Y))$ such that $q$ 
is a quotient space of $s$ and $dim_{\bf R}q=dim_{\bf Q_p}s$.
Let $\tilde Y$ be a completion of $Y$ and $Q_{\alpha }$
be polyhedra of the polyhedral expansion of $\tilde Y$.
Then there are uniform polyhedra $P_{\alpha }$ such that
${\tilde X}=\lim \{ P_{\alpha },
f^{\alpha }_{\beta }, E \} $ and
${\tilde \theta }({\tilde X})=\lim \{ Q_{\alpha },
{\tilde f}^{\alpha }_{\beta }, E \} ={\tilde Y},$
where ${\tilde f}^{\alpha }_{\beta }:={\tilde 
\theta }_{\beta }\circ f^{\alpha }_{\beta }\circ 
{\tilde \theta }_{\alpha }^{-1}:$ and then analogously 
to the proof of Theorem 5.1 we get statements of Theorem 5.2,
where $X={\tilde \theta }^{-1}(Y)$.
\par {\bf 5.3. Corollary.} {\it Let spaces $X$ and $Y$ be infinite
as in Theorems 5.1, 5.2, then for absolute polyhedral expansions
$T= \{ P_{\alpha }, f^{\alpha }_{\beta }, {\sf U} \} $ of $X$ and
$S= \{ Q_{\alpha }, g^{\alpha }_{\beta }, {\sf V} \} $ of $Y$
there is equality $card ({\sf U})=card ({\sf V})$.}
\par {\bf Proof.} In view of Theorems 5.1 and 5.2
to the inverse mapping system 
$T_X= \{ P_{\alpha }, f^{\alpha }_{\beta }, {\sf U}_X \} $ 
there corresponds the inverse mapping system
$T_Y= \{ Q_{\alpha }, {\tilde f}^{\alpha }_{\beta }, {\sf U}_X \} $ 
such that $X=\lim T_X$ and $Y=\lim T_Y$, where $P_{\alpha }$ and
$Q_{\alpha }$ are polyhedra over $\bf L$ and $\bf R$ respectively. 
The bonding mappings
$f^{\alpha }_{\beta }$ are simplexial, but ${\tilde f}^{\alpha }_{\beta }$
are not simplexial. Using permissible modifications 
(see \S 2.8, Lemmas 3.2 and 3.8) we get a polyhedral expansion
$W = \{ Z_{\alpha }, h^{\alpha }_{\beta }, F \} $ of $Y$
such that $h^{\alpha }_{\beta }$ are simplexial and $Z_{\alpha }$ 
are polyhedra. Moreover, $card (F)=card ({\sf U}_X)\aleph _0=
card ({\sf U}_X)$, since polyhedra $Q_{\alpha }$ are metrizable.
Using sections of the polyhedral expansion for $Y$ 
we get $S$ with $card({\sf U})=card({\sf V})$.
Vice versa, using $S$ we get $S_X$, then with the help of permissible
modifications we construct by $S_X$ polyhedral expansion of $X$ and 
using sections of spectrum we get an absolute polyhedral expansion for $X$.
\par {\bf 5.4. Note.} Certainly there is not 
in general a mapping of $T$ into $S$ or vice versa for $T$ and $S$ from
\S 5.3.
\par Each simplex $s$ is clopen in its polyhedra $P$
over $\bf L$, but in general $P$ may be not open in the corresponding 
Banach space, because may be cases of $s$ in $P$ 
such that $dim_{\bf L}s<dim_{\bf L}P=\sup_{s\subset P}dim_{\bf L}s$. 
There is the following particular case of spaces $X$,
when in their polyhedral expansions $S= \{ P_n, f^n_m, E \} $
polyhedra $P_n$ are clopen in the corresponding Banach spaces.
\par {\bf 5.5.} Let $H$ be a locally $\bf K$-convex space, 
where $\bf K$ is a non-Archimedean field. Let $M$ be a manifold modelled 
on $H$ and $At(M) = \{ (U_j,f_j):$ $ j\in A \} $ be an atlas
of $M$ such that $card (A)\le w(H)$, where $f_j: U_j\to V_j$
are homeomorphisms, $U_j$ are open in $M$, $V_j$ are open in $H$,
$\bigcup_{j\in A}U_j=M$, $f_i\circ f_j^{-1}$ are continuous on 
$f_j(U_i\cap U_j)$ for each $U_i\cap U_j\ne \emptyset $.
Let $\bf \tilde K$, $\tilde H$ and $\tilde M$ denote completions
of $\bf K$, $H$ and $M$ relative to their uniformities.
\par {\bf Theorem.} {\it If $H$ is infinite-dimensional over $\bf K$
or $\bf \tilde K$ is not locally compact, 
then $M$ is homeomorphic to the clopen subset of $H$.}
\par {\bf Proof.} Since $\tilde H$ is the complete
locally $\bf \tilde K$-convex space, then
${\tilde H}=pr-\lim \{ {\tilde H}_q, \pi ^q_v, F \} $ is a projective limit
of Banach spaces ${\tilde H}_q$ over $\bf \tilde K$, where $q\in F$, 
$F$ is an ordered
set, $\pi ^q_v: {\tilde H}_q\to {\tilde H}_v$ 
are linear continuous epimorphisms.
Therefore, each clopen subset $W$ in $\tilde H$ has a decomposition
$W=\lim \{ W_q, \pi ^q_v, F \} $, where $W_q=\pi ^q_v(W)$ are clopen in
${\tilde H}_q$. The base of topology of $\tilde M$ 
consists of clopen subsets.
If $W\subset V_j$, then $f_j^{-1}(W)$ has an analogous decomposition.
From this and Proposition 2.5.6 \cite{eng} it follows, that
$\tilde M=\lim \{ {\tilde M}_q, {\tilde \pi }^q_v, F \} $, where 
${\tilde M}_q$ are manifolds
on ${\tilde H}_q$ with continuous bonding mappings between charts
of their atlases.
If $H$ is infinite-dimensional over $\bf K$, then each ${\tilde H}_q$ is 
infinite-dimensional over $\bf \tilde K$ \cite{nari}.
From $card(A)\le w(H)$ it follows, that each ${\tilde M}_q$ has an atlas
$At'({\tilde M}_q)=\{ {U'}_{j,q};f_{j,q};{A'}_q \} $ 
equivalent to $At({\tilde M}_q)$ such that
$card ({A'}_q)\le w(H_q)=w({\tilde H}_q)$, since $w({\tilde H})=w(H)$,
where $At({\tilde M}_q)$
is induced by $At({\tilde M})$ by the quotient mapping
${\tilde \pi }_q: {\tilde M}\to {\tilde M}_q$.
In view of Theorem 2 \cite{lu4} each ${\tilde M}_q$ is homeomorphic 
to a clopen subset ${\tilde S}_q$ of ${\tilde H}_q$, 
where $h_q: {\tilde M_q}\to {\tilde S_q}$ are homeomorphisms. 
To each clopen ball $\tilde B$ in ${\tilde H}_q$ there corresponds
a clopen ball $B={\tilde B}\cap H_q$ in $H_q$, hence
$S_q={\tilde S}_q\cap H_q$ is clopen in $H_q$ and
$h_q: M_q\to S_q$ is a homeomorphism.
Therefore, $M$ is homeomorphic to
a closed subset $V$ of $H$, where $h: M\to V$ is a homeomorphism,
$V\subset H$, $h=\lim \{ id, h_q, F \} $, $id: F\to F$
is the identity mapping. Since each $h_q$ is surjective, then
$h$ is surjective by Lemma 2.5.9 \cite{eng}.
If $x\in M$, then ${\tilde \pi }_q(x)=x_q\in M_q$,
where $\pi : H\to H_q$ are linear quotient mappings
and ${\tilde \pi }_q: M\to M_q$ are induced quotient mappings.
Therefore, each $x\in M$ has a neighbourhood ${\tilde \pi }_q^{-1}(
Y_q)$, where $Y_q$ is an open neighbourhood of $x_q$ in $M_q$.
Therefore, $h(M)=V$ is open in $H$.
\par {\bf 5.6. Theorem.} {\it Let $M$ be a manifold satisfying 
conditions of Theorem 5.5, where $\bf K$ is with non-trivial
valuation, $M$ and $\bf K$ are complete relative to their own uniformities. 
Then $M$ has a polyhedral expansion
$S=\{ P_n, f^n_m, E \} $ consisting of polyhedra $P_n$ clopen in 
the corresponding Banach spaces $H_q$ over $\bf K$.}
\par {\bf Proof.} We treat the case of non-complete
$\bf K$ also.
Let $\bf \tilde K$ be a completion of $\bf K$
relative to its uniformity. 
If $\bf \tilde K$ is not locally compact
then $\bf \tilde K$ contains a proper locally compact 
subfield $\bf \tilde L$ \cite{roo}. Since $M$ is homeomorphic to
a clopen subset of $H$, then a completion $\tilde M$ of $M$
is homeomorphic to a clopen subset in a completion
$\tilde H$ of $H$, where $\tilde H$ is over $\bf \tilde K$. 
If $\bf \tilde L$ is a locally compact subfield of $\bf \tilde K$,
then $\bf \tilde L$ is spherically complete and $\tilde H$
can be considered as the topological vector space
over $\bf \tilde L$ \cite{nari,roo}.
By Lemmas and Theorems of \S 3 we can choose $P_n$ clopen 
in the corresponding Banach spaces over $\bf \tilde L$.
We also can choose polyhedra over $\bf \tilde K$,
repeating proofs for $\tilde M$ over $\bf \tilde K$.
In \S 3 polyhedra over locally compact field were chosen
that to encompass cases of compact spaces $X$, but 
polyhedra can be considered also over $\bf \tilde K$.
For compact $X$ and non locally compact $\bf \tilde K$
such polyhedral expansions have little sense, since 
discrete $P_n$ with singletons $s$ as simplexes may appear
such that $dim_{\bf \tilde K}s=0$.
But for the manifold $\tilde M$ due to embedding into
$\tilde H$ as the clopen subset we get the polyhedral expansion
such that $P_n$ are non-trivial over $\bf \tilde K$, moreover
there exists a polyhedral expansion such that $P_n$ are clopen 
in the corresponding Banach spaces ${\tilde H}_q$
for each $n$.
Using absolute polyhedral expansions constructed with the help of 
sections of the initial polyhedral expansion of $\tilde M$ from \S 4, 
we get  polyhedral expansions
of $M$ over $\bf \tilde L$ and also over $\bf \tilde K$.
In a particular case of complete $\bf L$ and $\bf K$ relative to their
uniformities the manifold $M$ is complete together with $H$, that is
$\tilde M=M$ and $\tilde H=H$, consequently, $M$ as the topological space 
is homeomorphic to $\lim S$.
\par The author expresses his gratefulness to Professor G. Tironi
for his interest in this job and hospitality at Mathematical Department
of Trieste University since 24 May until 24 October 1999,
also to Professor V.V. Fedorchuk at the Mathematical Department
of the Moscow State University for careful reading of the manuscript.


\newpage
\begin{thebibliography}{99}
\bibitem{coris} H.H. Corson, J.R. Isbell. "Some properties of 
strong uniformities." Quart. J. Math. 1960, V. {\bf 11}, 17-33.
\bibitem{steil} S. Eilenberg, N. Steenrod. "Foundations of 
algebraic topology" (Princ. Univ. Press, Princeton, New Jersey, 1952).
\bibitem{eng} R. Engelking. "General topology" (Mir: Moscow, 1986).
\bibitem{fed} V.V. Fedorchuk, A.Ch. Chigogidze. "Absolute retracts
and infinte-dimensional manifolds" (Nauka: Moscow, 1992).
\bibitem{freu} H. Freudenthal. "Entwicklungen von R\"aumen und
ihren Gruppen." Compositio Mathem. 1937, V. {\bf 4}, 145-234.
\bibitem{isb1} J.R. Isbell. "Euclidean and weak uniformities."
Pacif. J. Math. 1958, V. {\bf 8}, 67-86.
\bibitem{isb2} J.R. Isbell. "On finite-dimensional uniform spaces."
Pacif. J. Math. 1959, V. {\bf 9}, 107-121.
\bibitem{isb3} J.R. Isbell. "Irreducible polyhedral expansions."
Indagationes Mathematicae. Ser. A. 1961, V. {\bf 23}, 242-248.
\bibitem{isb4} J.R. Isbell. "Uniform neighborhood retracts."
Pacif. J. Math. 1961, V. {\bf 11}, 609-648.
\bibitem{isb5} J.R. Isbell. "Uniform spaces." Mathem. Surveys.
1964, V. {\bf 12}, AMS. Providence, R.I., USA.
\bibitem{itrsw} G. Itzkowitz, S. Rothman, H. Strassberg, T.S. Wu.
"Characterisation of equivalent uniformities in topological groups."
Topology and its applications. 1992, V. {\bf 47}, 9-34.
\bibitem{koz} I.M. Kozlovsky. "Absolute polyhedral expansions
of metric spaces." Trudy Moskov. Mat. Obscestva. 1979, V. {\bf 40},
83-119.
\bibitem{lu1} S.V. Ludkovsky. "Free locally convex spaces generated by
locally compact groups." Usp. Mat. Nauk. 1993, V. {\bf 48}, N 2, 189-190.
\bibitem{lu2} S.V. Ludkovsky. "Free locally convex spaces and their
isomoprphisms." Usp. Mat. Nauk. 1994, V. {\bf 49}, N 6, 207-208.
\bibitem{lu3} S.V. Ludkovsky. "Non-Archimedean free Banach spaces."
Fundam. i Prikl. Mathem. 1995, V. {\bf 1}, 979-987.
\bibitem{lu4} S.V. Ludkovsky. "Embeddings of non-Archimedean Banach 
manifolds into non-Archimedean Banach spaces."
Usp. Mat. Nauk. 1998, V. {\bf 53}, N 5, 241-242.
\bibitem{lu5} S.V. Ludkovsky. "Non-Archimedean polyhedral 
expansions of ultrauniform spaces."
Usp. Mat. Nauk. 1999, V. {\bf 54}, N 5, 163-164.
\bibitem{monzi} D. Montgomery, L. Zippin. "Topological transformation 
groups" (J. Wiley and Sons: New York, 1955).
\bibitem{nari} L. Narici, E. Beckenstein. "Topological vector spaces"
(Marcel-Dekker Inc., New York, 1985).
\bibitem{pas} B.A. Pasynkov. "About spectral expansions of 
topological spaces." Mathem. Sborn. 1965, V. {\bf 66}, N 1, 35-79. 
\bibitem{pont} L.S. Pontryagin. "Continuous groups" (Nauka: Moscow, 1984).
\bibitem{roo} A.C.M. van Rooij. "Non-Archimedean functional analysis"
(Marcel Dekker: New York, 1978).
\bibitem{wei} A. Weil. "Basic number theory" Sec. Edit. 
(Springer: Berlin, 1973).
\end{thebibliography}
\end{document}